\documentclass[12pt]{article}
\usepackage{amsmath}
\usepackage{graphicx}
\usepackage{natbib}
\usepackage{url} % not crucial - just used below for the URL 

\usepackage[colorlinks=true,linkcolor=blue]{hyperref}%

\hypersetup{urlcolor=blue, citecolor=blue, colorlinks=true, linkcolor=blue} 

\usepackage{multicol}
\usepackage{bm}

\usepackage{amssymb}
\usepackage{amsfonts}

\usepackage{tikz}
\usetikzlibrary{decorations.pathreplacing}

\usepackage{authblk}

\newtheorem{proposition}{Proposition}[section]

\newtheorem{prueba}{Proof}[section]

\newtheorem{remark}{Remark}[section]

\usepackage{algorithm}
\usepackage{algorithmicx}
\usepackage{algpseudocode}
\begin{document}

\def\spacingset#1{\renewcommand{\baselinestretch}%
{#1}\small\normalsize} \spacingset{1}

%%%%%%%%%%%%%%%%%%%%%%%%%%%%%%%%%%%%%%%%%%%%%%%%%%%%%%%%%%%%%%%%%%%%%%%%%%%%%%

\title{\bf Computationally Efficient Algorithms for Simulating  Isotropic Gaussian Random Fields on Graphs with Euclidean Edges}

\author[1]{Alfredo Alegría}
\affil[1]{Departamento de Matemática, Universidad Técnica Federico Santa María, Valparaiso, Chile}

\author[2,3]{Xavier Emery}
\affil[2]{Department of Mining Engineering\\ Universidad de Chile\\ Santiago, Chile}
\affil[3]{Advanced Mining Technology Center\\ Universidad de Chile\\ Santiago, Chile}

\author[4]{Tobia Filosi}
\affil[4]{Department of Mathematics, University of Trento, Trento, Italy} 

\author[5,6]{Emilio Porcu}
\affil[5]{Department of Mathematics, and BTC Center for Biotechnology, Khalifa University }
\affil[6]{ADIA LAB. Abu Dhabi, United Arab Emirates
}

\maketitle

%\if0\blind
%{
%  \title{\bf Computationally Efficient Algorithms for Simulating  Isotropic Gaussian Random Fields on Graphs with Euclidean Edges}
%  \author{Alfredo Alegr\'ia\\
%    Departamento de Matem\'atica, Universidad T\'ecnica Federico Santa Mar\'ia,\\ Valparaiso, Chile,\\
%     Xavier Emery\\  Department of Mining Engineering, Universidad de Chile, Santiago, Chile,\\
%  Advanced Mining Technology Center, Universidad de Chile, Santiago, Chile,  \\
%     Tobia Filosi\\ 
%     Department of Mathematics, 
%     University of Trento, Trento, Italy\\
%  and\\
%   Emilio Porcu\\
%   Department of Mathematics, and BTC Center for Biotechnology, Khalifa University\\  
%ADIA LAB. Abu Dhabi, United Arab Emirates
%}
    
%  \maketitle
%} \fi

%\if1\blind
%{
%  \bigskip
%  \bigskip
%  \bigskip
%  \begin{center}
%    {\LARGE\bf Title}
%\end{center}
%  \medskip
%} \fi

\bigskip
\begin{abstract}
 This work  addresses the problem of simulating Gaussian random fields that are continuously indexed over a class of metric graphs, termed graphs with Euclidean edges, being more general and flexible than linear networks. We introduce three general algorithms that allow to reconstruct a wide spectrum of random fields having a covariance function that depends on a specific metric, called resistance metric, and proposed in recent literature. 
    The algorithms are applied to a synthetic case study consisting of a street network. They prove to be fast and accurate in that they reproduce the target covariance function and provide random fields whose finite-dimensional distributions are approximately Gaussian.
\end{abstract}

\noindent%
\vfill

\newpage
\spacingset{1.75} % DON'T change the spacing!

\section{Introduction}

Random fields having networks as an index set are ubiquitous across various applications, such as assessing ecological and environmental impacts in road networks, understanding the propagation of signals through dendrites, or modeling the flow of rivers, communication networks and electrical wires. For a recent comprehensive review on this subject, readers are directed to \cite{baddeley2021analysing}.  Several applications involving non-standard networks are discussed by \cite{porcu2023stationary}.

The statistical analysis of network data poses significant challenges. Utilizing the coordinates of points on the network akin to coordinates in a Euclidean space where the network is embedded is conceptually flawed, distorting the proximity among points, and potentially producing unrealistic forecasts at locations of the network that have not been observed.

Our paper deals with the simulation of spatial random fields having a (not necessarily linear) network as an index set. In our pursuit of formulating simulation methodologies, our aim is to integrate the extensive body of knowledge and practicality offered by linear networks, while also addressing the challenges posed by more intricate network structures. Consequently, the findings presented throughout this manuscript retain their validity across various network architectures, as described below, thereby expanding their score of applicability. This general formulation can prove useful in situations where  linear network structures may not accurately represent the complexity of the spatial domain \citep{ver2006spatial, ver2010moving, baddeley2021analysing, porcu2023stationary}.

The problem of simulation {\em on} networks is apparently different from that of simulation {\em of} networks. For the latter, there is a massive literature at hand, coming from computer science, machine learning and applied mathematics. For the former, the literature is sparse as outlined below.

For a Gaussian random field defined over any network, it is customary to assume that the covariance function is isotropic \citep{anderes2020isotropic}, that is, the covariance between any pair of observations located over two different points of the network is solely a function of the \emph{distance} between the points. While this concept is intuitive and natural for random fields defined over Euclidean spaces or manifold, defining distances over networks is tricky. The recent works of \cite{anderes2020isotropic}, \cite{bolin2024gaussian}, \cite{filosi2023temporally} and \cite{porcu2023stationary} illustrate the intricacies behind the implementation of the right metrics for a random field defined over a non-linear network. The problem can be understood by noting that, essentially, a non-linear network is a (semi) metric space, that is a non-empty set endowed with a (semi) metric. 

This paper challenges the problem of simulation of (approximately) Gaussian random fields having a metric graph as an index set. Specifically, the type of metric graph considered in this paper is termed a {\em graph with Euclidean edges} after the tour de force by \cite{anderes2020isotropic}. These graphs allow for a continuously defined random field where points are defined in both vertices and arcs. 
This framework proves valuable for analyzing domains featuring tunnels or bridges, incorporating weighted edges within the graph structure in a cohesive manner.

\cite{anderes2020isotropic} prove that Gaussian random fields on graphs with Euclidean edges can be endowed with, either, the geodesic distance --- provided some technical constraints are fulfilled --- or the resistance metric, being an extension of the conventional metric used within the framework of electrical networks \citep{Klein}. The latter metric not only provides a physical interpretation but also allows for the establishment of more flexible correlation structures on the graph, accommodating a broader range of real-world applications, unlike the more restrictive geodesic distance.

The simulation of random fields over non-linear networks has been challenged to a limited extent only. \cite{moller2022cox} propose simulation algorithms on linear networks that employ operations of large covariance matrices, such as the Cholesky decomposition, and are challenging from a computational viewpoint when the sample size is large. \cite{bolin2024gaussian} studied a spectral decomposition for random fields on metric graphs. This approach depends on the eigenfunctions and eigenvalues of a differential operator, covering the case of a covariance function analogous to the Mat{\'e}rn model and depending on geodesic distance. For  other families of covariance functions, this approach is challenging, if not unfeasible. \cite{bolin2023regularity} provide sampling techniques to generalized Whittle-Mat\'ern fields  with general smoothness based on finite element approximations and rational approximations. The method is fast as it uses sparse matrix methods, and it is essentially linear in complexity for most graphs. It turns to be even faster than sampling using the stochastic partial differential equation (SPDE) approach in 2D.  \cite{bolin2023statistical} provide simulation techniques for Whittle-Mat\'ern fields with exact integer smoothness. Using the Markov property of the field allows for fast simulation. 

Although the approaches provided above offer insightful techniques that are computationally fast and reliable, we note that they are limited to random fields associated with specific classes of SPDE. This is not an issue {\em per se}. Yet, alternative paths might allow for exploring simulations of  random fields that are not necessarily associated with specific SPDEs and their Markov approximations.

While the marriage between SPDEs and Gaussian Markov random fields (their discrete counterpart) is substantially based on the precision matrix (inverse of the covariance matrix), this paper explores different paths. We exploit the connection between resistance metric with the variogram of certain classes of {\em ad hoc}  random fields in concert with three alternative stochastic representations. 
This leads to three efficient simulation algorithms 
inspired by classical geostatistical approaches that have been originally proposed to simulate isotropic random fields (with respect to Euclidean distance) in $\mathbb{R}^d$, for $d$ a positive integer. Our proposal is rich: we are able to simulate Gaussian random fields from a wealth of choices of covariance functions depending on the resistance metric.

We then proceed to illustrate algorithmic effectiveness through examples built upon existing popular networks found in R packages, such as the \texttt{spatstat} package \citep{baddeley2005spatstat}. Our simulation routines are compatible with the syntax and structure of networks within these packages. The codes are accessible in the GitHub repository \url{https://github.com/alfredoalegria/FastSimNetworks} and the results within the manuscript are reproducible.

The paper is structured as follows. Section \ref{preliminares} presents preliminaries on graphs with Euclidean edges and the resistance metric. Section \ref{constructions} introduces three stochastic representations for random fields on the graph, serving as the foundation for the simulation algorithms described in Section \ref{simulation}. Section \ref{ejemplos} illustrates several numerical examples and offers guidelines for practitioners. Finally, Section \ref{discussion} concludes the paper with a discussion.

\section{Preliminaries}
\label{preliminares}

\subsection{Graphs with Euclidean Edges}

The following lines provide a succinct and non-rigorous illustration of graphs with Euclidean edges. For a formal treatment, the reader is referred to \cite{anderes2020isotropic} and \cite{porcu2023stationary}. The recent works related to random fields that are continuously defined over networks are substantially based on the following facts:  
\begin{enumerate}
    \item A non-linear network can be represented through a special graph structure as depicted below. The graph allows for non-standard topologies that elude the limits of linear network representations.
    \item In turn, the graphs proposed in the recent literature, endowed with a proper (semi) metric, become (semi) metric spaces.
    \item (Semi) metric spaces can be {\em embedded} into convenient function spaces where implementing covariance functions become straightforward thanks to the contributions from the early 1940ies \citep{schoenberg}.
\end{enumerate}
The type of topology considered by \cite{anderes2020isotropic} is termed that of a graph {\em with Euclidean edges}, denoted with a triple $\mathcal{G} = \left( \mathcal{V},\mathcal{E},\{\varphi_e\}_{e\in\mathcal{E}} \right)$ where the elements are blended in the following way: 
\begin{enumerate}
    \item[(a)] $(\mathcal{V},\mathcal{E})$ has a graph structure, where $\mathcal{V}$ is the set of vertices and $\mathcal{E}$ contains the edges. We assume that this graph is simple and connected, {i.e.}, $\mathcal{V}$ is finite, the graph has not repeated edges or edges that join a vertex to itself, and every pair of vertices is connected by a path.  
    \item[(b)] Each edge $e\in\mathcal{E}$ is associated with a unique abstract set, also denoted $e$, such that $\mathcal{V}$ and all the edge sets are mutually disjoint. 
    \item[(c)] Let  $u$ and $v$ be vertices connected by $e\in \mathcal{E}$. Then, $\varphi_e$ is a bijective mapping defined on $e \cup \{u,v\}$, such that $\varphi_e$ maps $e$ onto an open interval $(\underline{e},\overline{e}) \subset \mathbb{R}$ and $\{u,v\}$ onto $\{\underline{e},\overline{e}\}$.
\end{enumerate}

Throughout, when we write $u\in\mathcal{G}$, we are referring to an element that can be either a vertex or a point within an edge set, i.e., $u\in \mathcal{V}\cup \bigcup_{e\in\mathcal{E}} e$.

Some comments are in order. A linear network is formally defined as the union of linear segments, denoted as $\ell_i \subset \mathbb{R}^2$, $i\in I$, and these segments intersect only at their endpoints. Apparently, the definition above shows that a linear network is a special case of a graph with Euclidean edges: the endpoints of the segments constitute the set of vertices, while each edge $e_i \in \mathcal{E}$ is defined as the relative interior of the linear segment $\ell_i$. The bijection $\varphi_{e_i}$ corresponds to the inverse of the length-path parameterization of $\ell_i$. However, graphs with Euclidean edges allow for far more flexibility: they can effectively model, for instance, bridges or tunnels. To see some examples of graphs with Euclidean edges that are not linear networks, refer to \cite[Figure 3]{anderes2020isotropic}. Insightful graphical representations of these graphs can also be found in \cite{filosi2023temporally}.

An additional recall to notation will provide an easier exposition throughout. Two vertices $u,v\in\mathcal{V}$ are neighbors  
if there exists an edge in $\mathcal{E}$ connecting them. If $u \in e\in\mathcal{E}$, we let $\underline{u}$ and $\overline{u}$ represent the vertices that are connected by $e$, such that $\underline{u}$ and $\overline{u}$ correspond, respectively, to $\underline{e}$ and $\overline{e}$. When $u\in\mathcal{V}$, we define $\underline{u}=\overline{u} = u$.

\subsection{Random Fields on Graphs}

The definition of random fields is broad and entails a wealth of specific cases. Here, random field stands for a stochastic process having a continuous index set. Specifically, a random field over ${\cal G}$ is an uncountable collection $Y:= \{Y(u), \; u \in {\cal G} \}$ of random variables in the same probability space. A random field $Y$ is called Gaussian when the finite-dimensional distributions of $(Y(u_1),\ldots,Y(u_n))^{\top}$, with $\top$ denoting the transpose of a vector and $n$ a positive integer, are multivariate Gaussian, for all $u_1,\hdots,u_n\in\mathcal{G}$, with a given mean vector and a covariance matrix, denoted $\Sigma= \left ( \Sigma_{ij} \right )_{i,j=1}^n$ throughout, with $\Sigma_{ij}= {\rm cov} \left ( Y(u_i),Y(u_j) \right )$ and {\rm cov} denoting {\em covariance}. 

For Gaussian random fields, modeling, inference and prediction are solely determined by the mean and the covariance function $K: {\cal G} \times {\cal G} \to \mathbb R$, which is a positive semidefinite mapping. The link between the function $K$ and the matrix $\Sigma$ is that $\Sigma_{ij} = K(u_i,u_j)$, for $u_i,u_j \in {\cal G}$. The function $K$ is called isotropic for the graph ${\cal G}$ when there exists a pair $(C, d_\cdot)$, with $d_\cdot$ a suitable semi metric, such that $K(u,v) = C(d_\cdot(u,v))$, and for a mapping $C: \sigma_{d_{\cdot}} \to \mathbb R$, with $\sigma_{d_{\cdot}}$ being the disk of $d_{\cdot}$, that is the image space for the mapping $d_\cdot$. 

Apparently, once the suitable (semi) metric $d_\cdot$ is found, the pair $({\cal G},d_\cdot)$ becomes a (semi) metric space. As for the functions $C$ such that $C(d_\cdot(\cdot,\cdot))$ is positive semidefinite, we defer this discussion as we first define the proper metric for the graph used in this paper. 

\subsection{Endowing the Graph with the Resistance Metric}

While the geodesic (aka, shortest-path) distance is a physically natural candidate for a metric over a graph with Euclidean edges, the works of \cite{anderes2020isotropic} and \cite{porcu2023stationary} show that such a choice resents from a collection of drawbacks and structural limitations when implementing a suitable covariance function based on that metric. Further, in a number of cases, some graph topologies are not suitable to the geodesic distance, because some properties as illustrated above no longer hold. This fact motivated \cite{anderes2020isotropic} to search for a metric that is methodologically more flexible and computationally more efficient, the latter meaning the computational burden associated with the computation of the Laplacian matrix for the graph is alleviated when using such a distance. A solution is found through a generalization of the classical resistance metric that is reminiscent to electrical networks \citep{Klein}.

Formally, the resistance metric on $\mathcal{G}$ is defined as
\begin{equation}
    \label{metric}
    d_R(u,v) = \text{var}(Z_{\mathcal{G}}(u) - Z_{\mathcal{G}}(v)), \qquad u,v\in\mathcal{G}, 
\end{equation}
where $Z_\mathcal{G}$ is an auxiliary random field on $\mathcal{G}$, adopting the following form
\begin{equation}
\label{Z_G}
    Z_\mathcal{G}(u) =  Z_\mu(u) + \sum_{e\in \mathcal{E}} Z_e(u), \qquad u\in\mathcal{G}.
\end{equation}
Let us establish a precise definition of the random fields $Z_\mu$ and $Z_e$ involved in (\ref{Z_G}):
\begin{itemize}
    \item Firstly, the random field $Z_\mu$ is defined on the vertices $v_1,\hdots,v_n \in \mathcal{V}$. Indeed, we endow it with a Gaussian structure: $(Z_\mu(v_1),\hdots,Z_\mu(v_n))^\top \sim \mathcal{N}(0,L^{-1})$, where $L$, described next, is a variation of the Laplacian matrix. Consider the function $c:\mathcal{V}\times \mathcal{V}\rightarrow [0,\infty)$ such that
    $c(u,v) = 1/d_G(u,v)$ when $u$ and $v$ are neighbors, with $d_G(u,v)$ being the length of the edge joining these points, %$u\overset{\cal G}{\sim} v$, 
    and $c(u,v) = 0$ otherwise, and the function $c(u) = \sum_v c(u,v)$.  Thus, the entries of $L$ can be computed according to the following structure:
    $$ L(u,v) = \begin{cases}
        1 + c(u_0) & \text{ if } u = v = u_0\\
        c(u)  & \text{ if } u=v\neq u_0\\
        -c(u,v) & \text{ otherwise, }
    \end{cases} $$
    where $u_0$ is an arbitrary vertex. The resulting matrix is strictly positive definite as a consequence of the extra term $1$ in an entry of the main diagonal \citep{anderes2020isotropic}.  
    
    \item Secondly, $Z_\mu$ is extended to the edges by employing a linear interpolation of the form 
    \begin{equation}
    \label{interpolation}
        Z_\mu(u) = (1-d(u))Z_\mu(\underline{u}) + d(u) Z_\mu(\overline{u}),
    \end{equation} 
    where
    $$d(u) = \begin{cases} 
          d_{G}(u,\underline{u})/d_{G}(\underline{u},\overline{u}) & \text{ if } u\notin \mathcal{V}(\mathcal{G})\\
           0  &  \text{ otherwise. }
    \end{cases}$$

    \item Finally, $Z_e(u) = B_e(\varphi_e(u))$ for $u\in e$, and $Z_e(u) = 0$ otherwise, with $B_e$ standing for standard independent Brownian bridges, being independent of $Z_\mu$, such that $B_e(\underline{e}) = B_e(\overline{e}) = 0$. 
    
\end{itemize}

The resistance metric exhibits some interesting properties. As a metric, it adheres to fundamental criteria: it is non-negative and symmetric, $d_{R}(u,v)=0$ if and only if $u=v$, and it satisfies the triangular inequality. Moreover, it maintains invariance under certain operations, such as splitting an edge into two (with the addition of an extra vertex along that edge's path) or merging edges at a vertex with two incident edges (effectively removing that vertex) \citep{anderes2020isotropic}.

\subsection{Which Functions $C$ can be Coupled with the Metric $d_R$?}

\cite{anderes2020isotropic} prove that a sufficient condition for a continuous mapping $C: \sigma_{d_R} \to \mathbb{R}$ to ensure $C(d_{R}(\cdot,\cdot))$ to be positive semidefinite on ${\cal G} \times {\cal G}$ is that $C$ is the restriction of a function $\psi$ that is completely monotone on the real line. These functions have a long history that traces back to \cite{schoenberg}. They are the Laplace transforms of positive and finite measures (see, e.g., Theorem 1.4 in \citealp{schilling_bernstein_2012}). More specifically, $\psi:[0,\infty)\to (0,\infty)$ is completely monotone if, and only if, there exists a (unique) finite measure $\mu$ on $[0,\infty)$ such that
\begin{equation}
    \label{eq:bernsteinRepresentation}
    \psi(x)=\int_0^{\infty} e^{-xt} \, \text{d}\mu(t), \quad x\in [0,\infty).
\end{equation}

As a consequence, each completely monotone function $\psi$ is strictly positive on its domain. Notice that the condition $\psi(0)=1$ is equivalent to $\mu(\infty)=1$, that is $\mu$ is the cumulative distribution function of a random variable on $[0,\infty)$.

\section{Useful Stochastic Representations}
\label{constructions}

\subsection{First Construction (Spectral Representation)}
 
We study the class of random fields on $\mathcal{G}$ defined according to
\begin{equation}
\label{eq1}
    Y(u) := \sqrt{-2 \ln(V)} \cos( Z_{\mathcal{G}}(u) W +  \Lambda ), \qquad u\in\mathcal{G},
\end{equation}
where $Z_{\mathcal{G}}$ is the Gaussian random field introduced in (\ref{Z_G}) and involved in the definition of the resistance metric. Here, $V \sim \text{Unif}(0,1)$, $\Lambda \sim \text{Unif}(0,2\pi)$ and $W\sim F$, with $F$ standing for a probability measure on $\mathbb{R}$ with no atom at the origin (i.e., $W$ is almost surely different from $0$). We assume that $Z_\mathcal{G}(\cdot)$, $V$, $\Lambda$ and $W$ are independent.

This construction is inspired by spectral and substitution methods developed in the geostatistics literature to simulate random fields in Euclidean spaces \citep{lantu2002, allard2020}. The following proposition demonstrates that this approach enables the generation of isotropic random fields with a broad spectrum of correlation functions.

\begin{proposition}
\label{prop_spectral}
    The random field introduced in (\ref{eq1}) possesses a zero mean. Its covariance function, denoted by $C_Y$, is isotropic with respect to the resistance metric and is given by
    \begin{equation}
        \label{cov0_spectral}
        C_Y\left( d_R(u,v)  \right) 
        =  \int_{\mathbb{R}} \exp\left( - d_R(u,v) \omega^2 / 2\right) F(\text{d}\omega), \qquad u,v\in \mathcal{G}.
    \end{equation}
\end{proposition}

\begin{prueba}
    Note that (\ref{eq1}) can be written as
    \begin{equation}
    \label{rewr}
       Y(u) = \sqrt{-2 \ln(V)} \big[ \cos( Z_{\mathcal{G}}(u) W) \cos( \Lambda )  - \sin( Z_{\mathcal{G}}(u) W) \sin( \Lambda ) \big]. 
    \end{equation}
    Thus, by 
    considering the fact that $\mathbb{E}(\cos \Lambda) = \mathbb{E}(\sin \Lambda) = 0$ and that $\Lambda$ is independent of $V$, $W$ and $Z_{\cal G}(\cdot)$, it becomes evident that the random field has a zero mean. By also accounting for the fact that $\mathbb{E}(\cos\Lambda  \sin\Lambda) = 0$, $\mathbb{E}(\cos^2\Lambda) = \mathbb{E}(\sin^2\Lambda) = 1/2$ and $\mathbb{E}(-\ln(V))=1$, one obtains
    \begin{equation}
    \label{cov_spectral}
    \text{cov}(Y(u),Y(v)) = \mathbb{E}\left[
    \cos\left(W  \{Z_\mathcal{G}(u) - Z_\mathcal{G}(v) \} \right) \right], \qquad u,v\in\mathcal{G},
    \end{equation}
    where expectation is taken with respect to $W$ and $(Z_{\mathcal{G}}(u), Z_{\mathcal{G}}(v))^\top$. 
    Since $(Z_{\mathcal{G}}(u) - Z_{\mathcal{G}}(v)) / \sqrt{d_R(u,v)}$ follows a standard Gaussian distribution,  (\ref{cov_spectral}) can be written as
    \begin{equation}
        \label{cov2_spectral}
        \text{cov}(Y(u),Y(v))
            = 
        \frac{1}{\sqrt{2\pi}} \int_{\mathbb{R}} \left[ \int_{\mathbb{R}} \cos\left(\sqrt{d_R(u,v)}  \omega z \right) \exp\left( -z^2/2 \right)\text{d}z \right]   F(\text{d}\omega).
    \end{equation}
    We conclude that the random field $Y(\cdot)$ is isotropic with respect to the resistance metric, and the expression (\ref{cov0_spectral}) is readily obtained by solving the integral within the brackets in (\ref{cov2_spectral}) (which is related to the real part of the characteristic function of a normal random variable) \citep[formula 3.896.4]{Grad}.
    
    \hfill
  $\blacksquare$  
\end{prueba}

Proposition \ref{prop_spectral} shows that the covariance function $C_Y$ is a scale mixture of exponential covariances, therefore it is a completely monotone function. Reciprocally, a straightforward change of variable in (\ref{cov0_spectral})  shows that any completely monotone covariance function can be obtained with the proposed construction.
%\tobi{please someone check this proposition (I have run simulations and it works). It should generalise results in Table 3.1.}  {Seems OK to me.} \emi{Some small mistakes - see the brown color below.}
%\begin{proposition}
%\label{complmonot}
%    Let $\psi: [0,\infty)\to (0,1]$ be a completely monotone correlation function and let $\mu$ the unique measure such that (\ref{eq:bernsteinRepresentation}) holds. Let $W$ be the nonnegative random variable with cumulative distribution function $F(\omega):=1_{\omega \geq 0}\,\mu\left(\frac{\omega^2}{2}\right)$, $\omega\in\mathbb R$, where $\mu$ is the measure defined in (\ref{eq:bernsteinRepresentation}). Then, the process $Y$ \emi{as being} defined in Equation (\ref{eq1}) has covariance function
%    \begin{equation*}
%        C_Y(d_R(u,v)) = \psi(d_R(u,v)), \quad u, v \in {\cal G}.
%    \end{equation*}
%\end{proposition}
%\begin{prueba}
    %Let $d:=d_R(u,v)$: b
%    By Equation (\ref{cov0_spectral}),
%    \begin{align*}
%        C_Y\left( d_R(u,v)  \right) 
%        &=  \int_{\mathbb{R}} \exp\left( - d_R(u,v) \,\omega^2 / 2\right) \text{d}F(\omega) \\%=\int_0^{+\infty} \exp \left( -d_R(u,v)\frac{\omega^2}{2} \right) \text{d} \mu\left( \frac{\omega^2}{2} \right) \\
 %       &= \int_0^{\infty} \exp \left( -d_R(u,v) t \right) \text{d} \mu\left( t \right) = \psi(d_R(u,v)),
 %   \end{align*}
 %   where in the second-to-last step we have changed the variable by means of the monotonic transformation $t:=\frac{\omega^2}{2}$.
 %   $\qed$
%\end{prueba}
In Table \ref{tab:cov}, we provide a list of parametric families of  correlation functions on $\mathcal{G}$, for different choices of the spectral measure $F$. Some of these families (entries 1, 6 and 8) have been reported in \cite{anderes2020isotropic}. The last seven entries have been established by using formulae 4.2.16, 4.2.25, 4.2.26, 4.4.15, 4.5.3, 4.5.23 and 4.5.28 of \cite{erdelyi} to calculate (\ref{cov0_spectral}).

\begin{table}
    \centering
        \caption{Parametric correlation functions on $\mathcal{G}$ in terms of the resistance metric. For each model $C_Y$, we also provide the associated spectral measure $F$. $\Gamma$ is the gamma function, erf is the Gauss error function, $K_\tau$ is the modified Bessel function of the second kind of order $\tau$, $\tau>0$ is a shape parameter, and $a>0$ is a scale parameter. }
        
    \begin{tabular}{ll} \hline \hline
    $F(\text{d}\omega)$  & $C_Y$   \\ \hline 
    $\text{d}\delta_{\omega}(\{a\})$     &  $C_Y(d_R) = \exp( - a^2 d_R/2)$        \\     
    %$(2a\sqrt{\pi})^{-1}\exp(- \omega^2/(4a^2)) d\omega$  &  $C_Y(d_R) = (1+2 a^2 d_R )^{-1/2}$   \\
     {$(2a)^{-1} 1_{|\omega|<a} d\omega$} &   {$C_Y(d_R) =\sqrt{\pi} (2 a^2 d_R)^{-1/2} \, \text{erf}(a\sqrt{ d_R/2})$}  \\
      $(\pi a)^{-1}(1+\omega^2/a^2)^{-1}\text{d}\omega$   &   {$C_Y(d_R) = \exp( a^2d_R/2 ) [1-\text{erf}(-a\sqrt{d_R})]$}       \\ %$C_Y(d_R) = 2\exp( a^2d_R/2 ) \Phi(-a\sqrt{d_R})$
     {$a (\pi |\omega|)^{-1}(\omega^2-a^2)^{-1/2} 1_{|\omega|>a} d\omega$} &   {$C_Y(d_R) = 1-  \text{erf}(a \sqrt{d_R/2}) $}  \\ 
     {$(|\frac{\omega}{a}|^3 1_{|\omega|<a} + |\frac{\omega}{a}|(2-(\frac{\omega}{a})^2) 1_{a<|\omega|<2a}) \frac{d\omega}{a}$} &   {$C_Y(d_R) = (a^2 d_R/2)^{-2} (1-\exp(-a^2 d_R/2))^2$}  \\ 
     {$a^\tau |\omega|^{2\tau-1} \exp(-a\omega^2) d\omega / \Gamma(\tau)$} &   {$C_Y(d_R) = (2a)^\tau (2a+d_R)^{-\tau} $}  \\ 
     {$(a\Gamma(1/4))^{-1} \sqrt{2}\exp(-\omega^4/(4a^4)) d\omega$} &   {$C_Y(d_R) = a \sqrt{d_R} \exp(a^4 d_R^2/8) K_{1/4}(a^4 d_R^2/8)/ \Gamma(1/4)$}  \\ 
     {$a\omega^{-2} \exp(-a^2/(4\omega^2))/(2\sqrt{\pi}) d\omega$} &   {$C_Y(d_R) = \exp(-a\sqrt{ d_R/2})$}  \\ 
    \hline \hline
    \end{tabular}
    \label{tab:cov}
\end{table}

%\textcolor{brown}{EP1: El problema principal de esta construcci{\'o}n es que solo permite strictly positive covariance functions. Al mismo tiempo, creo que no podemos evitar esto, a menos que no vayamos a un network lineal (en ese caso, todo vale).}  \\

%\textcolor{brown}{EP2: I checked 4.2.22 in \cite{erdelyi} and the result is different. Can you elaborate a little how you got the entry number 3 in the Table?}  {Thanks. It was 4.2.25 and a I fixed the formula.}

Concerning the  random field distributions, one has the following result.

\begin{proposition}
\label{prop_spectral2}
   The random field defined in (\ref{eq1}) has a standard Gaussian univariate distribution. Furthermore, its bivariate distributions are isotropic with respect to the resistance metric. They can be written as mixtures of bivariate Gaussian distributions and have an isofactorial representation with the Hermite polynomials as the factors: for $u, v \in {\cal G}, u \neq v$, the probability density function of $(Y(u),Y(v))$ can be expanded as
   \begin{equation}
    \label{isof}
    \begin{split}
        g_{u,v}(y,y^\prime) &= g(y) g(y^\prime) \sum_{\nu=0}^\infty \mathbb{E}\left[\rho(u,v)^\nu\right] H_\nu(y) H_\nu(y^\prime),\quad y, y^\prime \in \mathbb{R},
    \end{split}
    \end{equation}
    where $\rho(u,v) := \cos(W Z \sqrt{d_R(u,v)})$, $Z$ is a standard Gaussian random variable independent of $W$, $H_\nu$ is the normalized Hermite polynomial of degree $\nu$ and $g$ is the standard Gaussian probability density function. 
    %\emi{has isofactorial been defined?}
\end{proposition}

\begin{prueba}
    We follow \cite{matheron1982} to determine the Fourier transforms of the univariate and bivariate distributions of $Y(\cdot)$. Owing to the Box-Muller method \citep{Box}, one can rewrite (\ref{rewr}) as
    \begin{equation*}
    \label{rewr2}
         Y(u) = Y_1 \cos( Z_{\mathcal{G}}(u) W) + Y_2 \sin( Z_{\mathcal{G}}(u) W), \quad u \in {\cal G},
    \end{equation*}    
    where $Y_1$ and $Y_2$ are standard Gaussian random variables that are independent and independent of $W$ and $Z_{\cal G}(\cdot)$. Accordingly, the distribution of $Y(u)$ (with $u \in {\cal G}$) conditioned to $W$ and $Z_{\cal G}(\cdot)$ has the following Fourier transform, where $\mathsf{i}$ stands for the imaginary unit and $\lambda$ for a real value:
    \begin{equation}
    \label{rewr3}
    \begin{split}
       \mathbb{E}\left[\mathsf{e}^{\mathsf{i} \lambda Y(u)} \mid W, Z_{\cal G}(\cdot)\right] &= \exp\left[-\frac{\lambda^2}{2}\left(\cos^2(Z_{\cal G}(u) W) + \sin^2(Z_{\cal G}(u) W)\right) \right] \\
       &= \exp\left[-\frac{\lambda^2}{2} \right], \quad u \in {\cal G},
    \end{split}
    \end{equation}
    which does not depend on $W$ or $Z_{\cal G}(\cdot)$. It is deduced that the non-conditional distribution of $Y(u)$ is standard Gaussian, as its Fourier transform is the same as the last expression in (\ref{rewr3}). As for the bivariate distributions, the distribution of $(Y(u),Y(v))$ (with $u, v \in {\cal G}$) conditionally on $W$ and $Z_{\cal G}(\cdot)$ has the following Fourier transform, where $\lambda$ and $\eta$ stand for real values:
    % ciao Emilio
    \begin{equation*}
        \label{rewr4}
        \begin{split}
           \mathbb{E}\left[\mathsf{e}^{\mathsf{i} \lambda Y(u) + \mathsf{i} \eta Y(v)} \mid W, Z_{\cal G}(\cdot)\right] &= \exp\left(-\frac{\lambda^2+\eta^2}{2}-\lambda \eta \cos(Z_{\cal G}(u) W - Z_{\cal G}(v) W) \right), \quad u,v \in {\cal G},
        \end{split}
    \end{equation*}
    where $Z := (Z_{\mathcal{G}}(u) - Z_{\mathcal{G}}(v)) / \sqrt{d_R(u,v)}$ follows a standard Gaussian distribution and is independent of $W$. Hence, the non-conditional Fourier transform is
    \begin{equation*}
        \label{rewr5}
        \begin{split}
           \mathbb{E}\left[\mathsf{e}^{\mathsf{i} \lambda Y(u) + \mathsf{i} \eta Y(v)}\right] &= \mathbb{E}\left[\exp\left(-\frac{\lambda^2+\eta^2}{2}-\lambda \eta \cos(W Z \sqrt{d_R(u,v)}) \right)\right] \\
           &= \frac{1}{\sqrt{2\pi}} \int_{\mathbb{R}} \int_{\mathbb{R}} \exp\left(-\frac{\lambda^2+\eta^2}{2}-\lambda \eta \cos(\omega z \sqrt{d_R(u,v)}) \right) \mathsf{e}^{-\frac{z^2}{2}} \text{d}z \, F(\text{d}\omega), \quad u,v \in {\cal G},
        \end{split}
    \end{equation*}
    with the double integral being convergent insofar as the first exponential in the integrand is upper bounded by $1$.
    
    It is seen that the Fourier transform depends on $u$ and $v$ only through their mutual distance $d_R(u,v)$, which establishes that the bivariate distributions of $Y(\cdot)$ are isotropic with respect to the resistance metric, and that it is a scale mixture of bivariate Gaussian Fourier transforms. Given that $W$ is almost surely different from $0$, the joint probability density function of $(Y(u),Y(v))$ exists if $u \neq v$ and can be written as:
    \begin{equation*}
        g_{u,v}(y,y^\prime) = \mathbb{E}\left[\frac{1}{2\pi \sqrt{1-\rho(u,v)^2}} \exp\left(-\frac{y^2 - 2\rho(u,v) y  y^\prime + y^{\prime 2}}{2(1-\rho(u,v)^2)} \right)\right], \quad y, y^\prime \in \mathbb{R},
    \end{equation*}    
    with $\rho(u,v) := \cos(W Z \sqrt{d_R(u,v)})$. 
    The proof concludes by noting that the bivariate Gaussian distributions have an isofactorial representation with the Hermite polynomials as the factors (see \cite{lancaster} or \citet[p. 412]{Chiles2012}), which leads to the announced expansion (\ref{isof}); the convergence of this expansion is guaranteed because $\rho(u,v)$ belongs to the open interval $(-1,1)$ almost surely \citep{matheron1976}.
    %\begin{equation}
    %\label{isof}
    %\begin{split}
        %g_{u,v}(y,y^\prime) &= g(y) g(y^\prime) \sum_{\nu=0}^\infty \mathbb{E}\left[\rho(u,v)^\nu\right] H_\nu(y) H_\nu(y^\prime),\quad y, y^\prime \in \mathbb{R},
    %\end{split}
    %\end{equation}
    %where $H_\nu$ is the normalized Hermite polynomial of degree $\nu$ and $g$ is the standard Gaussian probability density function. 
    Using formulae 1.320.5, 1.320.7 and 3.896.4 of \cite{Grad}, one has:
    \begin{equation*}
       \mathbb{E}\left[\rho(u,v)^\nu\right] = 
       \begin{cases}
           \frac{2}{2^\nu} \mathbb{E}\left[\sum_{\kappa=0}^{\frac{\nu-1}{2}} \frac{\nu!}{\kappa!(\nu-\kappa)!} \exp\left(\frac{(\nu-2\kappa)^2 W^2 d_R(u,v)}{2}\right)\right] \text{ if $\nu$ is odd}\\
           \frac{\nu!}{2^\nu (\nu/2)!^2} + \frac{2}{2^\nu} \mathbb{E}\left[\sum_{\kappa=0}^{\frac{\nu}{2}-1} \frac{\nu!}{\kappa!(\nu-\kappa)!} \exp\left(\frac{(\nu-2\kappa)^2 W^2 d_R(u,v)}{2}\right)\right] \text{ if $\nu$ is even,} 
       \end{cases}
    \end{equation*}
    so that a sufficient condition for (\ref{isof}) to be expressed analytically is that the moment generating function of $W^2$ is defined on $\mathbb{R}$ and known (e.g., with the first entry of Table \ref{tab:cov}).

       \hfill
  $\blacksquare$ 
\end{prueba}

\subsection{Second Construction (Poisson Dilution Model)}

Consider the class of random fields on $\mathcal{G}$  given by
\begin{equation}
\label{random_token}
    Y(u) :=  \sum_{x \in X} \epsilon(x) f\left(Z_{\mathcal{G}}(u)-x\right), \qquad u\in\mathcal{G},
\end{equation}
where $f$ is real-valued square-integrable function on $\mathbb{R}$, $X$ is a stationary Poisson point process in $\mathbb{R}$ with intensity $\theta=1$ independent of $Z_{\cal G}(\cdot)$, $\{\epsilon(x): x\in\mathbb{R} \}$ are mutually independent weights, such that $\mathbb{P}(\epsilon(x)=1) = \mathbb{P}(\epsilon(x)=-1) = 1/2$ for each $x\in\mathbb{R}$, and are independent of $X$ and $Z_{\cal G}(\cdot)$. Let $\psi_f$ be the transitive covariogram of $f$, defined as
\begin{equation}
\label{trcov}
  \psi_f(h) = \int_{\mathbb{R}} f(x+h)f(x)\text{d}x, \quad h\in\mathbb{R}.  
\end{equation}
This construction adapts the dilution method from classical geostatistics \citep{lantu2002}. The following proposition demonstrates that this approach facilitates the generation of isotropic random fields and establishes the range of correlation functions achievable through this method.

\begin{proposition}
\label{prop_dilution}
   The random field defined in (\ref{random_token}) has a zero mean and is isotropic with respect to the resistance metric. Its covariance function, denoted as $C_Y$, is given by
    \begin{equation}
    \label{cov0}
        C_Y\left( d_R(u,v)  \right) 
    =  \frac{1}{\sqrt{2\pi}}  \int_{\mathbb{R}} \psi_f\left( \sqrt{d_R(u,v)} z \right) \exp(-z^2/2)\text{d}z.
    \end{equation}
\end{proposition}

\begin{prueba}
The random field clearly possesses a zero mean, insofar as $\mathbb{E}(\epsilon(x))=0$ for any $x \in X$ and $\epsilon(x)$ is independent of $Z_{\mathcal{G}}$ and $X$. Let $I_k = [c_k,d_k]$ denote a  sequence of closed intervals in $\mathbb{R}$, such that $I_k \subset I_{k+1}$ and $\cup_{k\in\mathbb{N}} I_k = \mathbb{R}$. The length of $I_k$ is denoted by $|I_k|$. We first analyze (\ref{random_token}), by restricting our attention to $X \cap I_k$.
 We proceed to calculate the covariance function of the random field
 $$ Y_{k}(u) :=  \sum_{x \in X \cap I_k} \epsilon(x) f\left(Z_{\mathcal{G}}(u)-x\right), \qquad u\in\mathcal{G}, k\in\mathbb{N}. $$
 By conditioning on the number and location of points of the Poisson process and utilizing the independence of the weights, we obtain
\begin{eqnarray*}
    \label{cov_2nd}
\text{cov}(Y_k(u),Y_k(v))
    & = & \mathbb{E}\left(  \sum_{n=0}^\infty \frac{\exp(- |I_k|)  |I_k|^n}{n!} \frac{n}{|I_k|}\int_{I_k}   f\left(Z_{\mathcal{G}}(u)-x\right) f\left(Z_{\mathcal{G}}(v)-x\right)  \text{d}x \right)\\
    & = &  \mathbb{E}\left(  \sum_{n= {1}}^\infty \frac{\exp(- |I_k|)  |I_k|^n}{n!} \frac{n}{|I_k|}\int_{Z_{\mathcal{G}}(v) - d_k}^{Z_{\mathcal{G}}(v)-c_k}   f\left(x+ Z_{\mathcal{G}}(u)-Z_{\mathcal{G}}(v) \right) f(x)  \text{d}x \right)\\
    & = &  \mathbb{E}\left(  \int_{Z_{\mathcal{G}}(v) - d_k}^{Z_{\mathcal{G}}(v)-c_k}   f\left(x+ Z_{\mathcal{G}}(u)-Z_{\mathcal{G}}(v) \right) f(x)  \text{d}x \right). 
\end{eqnarray*}
As $k$ approaches infinity, the result expands to $\mathbb{R}$; therefore, we obtain
\begin{equation*}
    \text{cov}(Y(u),Y(v))  =  \begin{cases}
    \psi_f(0) \text{ if $u=v$}\\
    \mathbb{E}\left[    \psi_f\left( \sqrt{d_R(u,v)} Z \right) \right] \text{ otherwise}, \\
    \end{cases}
\end{equation*}
where $Z:= (Z_{\mathcal{G}}(u)-Z_{\mathcal{G}}(v))/\sqrt{d_R(u,v)}$, which follows a standard Gaussian distribution.

    \hfill
  $\blacksquare$ 
\end{prueba}

In Table \ref{tab:cov2}, we provide a list of parametric families of  {correlation} functions on $\mathcal{G}$, for different choices of the function $f$.  {The entries of the table have been established by using formulae 3.321.2, 3.321.3, 3.321.4, 3.322.1 and 6.511.13 of \cite{Grad} and formulae 2.19, 7.27 and 7.28 of \cite{Chiles2012}.}

\begin{table}
    \centering
        \caption{Parametric  {correlation} functions on $\mathcal{G}$ in terms of the resistance metric, for different choices of the dilution function $f:\mathbb{R}\rightarrow \mathbb{R}$.  For each model, $a>0$ is a  {scale} parameter.}
        
    \begin{tabular}{ll} \hline \hline
    $f(t)$  & $C_Y$  \\ \hline 
     
    $a^{-1/2} 1_{|t| \leq a/2}$  & $C_Y(d_R) =    \text{erf}\left(\frac{a}{a\sqrt{2 d_R}}  \right)   -  \frac{2 - 2 \exp(- (a^2/2)/\sqrt{d_R})}{\sqrt{2\pi}}   \sqrt{d_R} $    \\
    $\exp(-a^2 t^2) (2/\pi)^{1/4} \sqrt{a}$ & $C_Y(d_R) = (1 + a^2 d_R)^{-1/2} $ \\
     {$K_0(a|t|) \sqrt{2a}/\pi $} &  {$C_Y(d_R) = \exp(a^2 d_R/2) (1-\text{erf}(a \sqrt{d_R/2})) $} \\
    \hline \hline
    \end{tabular}
    \label{tab:cov2}
\end{table}

\begin{remark}
    Unlike the conventional Euclidean case, this construction does not offer correlation functions with a compact support. Essentially, we take an additional expectation with respect to a Gaussian probability function, which prevents the occurrence of compactly supported models. 
\end{remark}

\subsection{Third Construction (Dilution of a Random Germ)}
\label{thirdconstruction}

Consider the class of {random fields} on $\mathcal{G}$  given by
\begin{equation}
\label{random_token2}
    Y(u) :=  \frac{\epsilon f\left(Z_{\mathcal{G}}(u)-X\right)}{\sqrt{\varpi(X)}}, \qquad u\in\mathcal{G},
\end{equation}
where $\epsilon$ is a random weight such that $\mathbb{P}(\epsilon=1) = \mathbb{P}(\epsilon=-1) = 1/2$, $X$ is a random variable with a probability density function $\varpi$ that is positive on $\mathbb{R}$, $f$ is real-valued square-integrable function on $\mathbb{R}$, and $(Z_{\cal G}(\cdot),\epsilon,X)$ are mutually independent. Let $\psi_f$ be the transitive covariogram of $f$, as per (\ref{trcov}).

\begin{proposition}
\label{prop_importancesampling}
   The {random field} defined in (\ref{random_token2}) has a zero mean, is isotropic with respect to the resistance metric, and its covariance function, denoted as $C_Y$, is given by (\ref{cov0}).
   % \begin{equation}
   % \label{cov0}
   %     C_Y\left( d_R(u,v)  \right) 
   % =  \frac{1}{\sqrt{2\pi}}  \int_{\mathbb{R}} \psi_f\left( \sqrt{d_R(u,v)} z \right) \exp(-z^2/2)\text{d}z.
   % \end{equation}
\end{proposition}

\begin{prueba}
The zero mean of $Y(u)$ for all $u \in {\cal G}$ stems from the zero mean of $\epsilon$ and the fact that $\epsilon$ is independent of $X$ and $Z_{\cal G}(\cdot)$. As for the covariance, one has:
\begin{equation*}
\begin{split}
\text{cov}(Y(u),Y(v))
    & =  \mathbb{E}\left( \int_{\mathbb{R}}   f\left(Z_{\mathcal{G}}(u)-x\right) f\left(Z_{\mathcal{G}}(v)-x\right)  \text{d}x \right)\\
    & =   \mathbb{E}\left(  \int_{\mathbb{R}}   f\left(x+ Z_{\mathcal{G}}(u)-Z_{\mathcal{G}}(v) \right) f(x)  \text{d}x \right)\\
    & =  \begin{cases}
     \psi_f(0) \text{ if $u=v$}\\   
    \mathbb{E}\left[    \psi_f\left( \sqrt{d_R(u,v)} Z \right) \right] \text{ otherwise},
    \end{cases}
\end{split}
\end{equation*}
where $Z:= (Z_{\mathcal{G}}(u)-Z_{\mathcal{G}}(v))/\sqrt{d_R(u,v)}$ follows a standard Gaussian distribution.

    \hfill
  $\blacksquare$ 
\end{prueba}

\begin{remark}
The class of covariance functions attainable with this construction is the same as the class obtained by the Poisson dilution model. With this version, however, we anticipate a simplified algorithm, as it no longer necessitates simulating a Poisson point process, which may involve aggregating many points and render its implementation less efficient.
\end{remark}

\section{Simulation Algorithms}
\label{simulation}

We exploit the stochastic representations of the previous section to propose specialized simulation algorithms. 

\subsection{Simulation of $Z_{\mathcal{G}}$}

The auxiliary {random field} $Z_{\mathcal{G}}$ is a key mathematical ingredient in the three constructions presented in the previous section. We employ the following steps to simulate this field:

\begin{itemize}

        \item[] {\bf Step 1.}  We first simulate $Z_\mathcal{G}$ on the vertices of the graph $v_1,\hdots,v_n$. It consists of  simulating a Gaussian random vector $Z_\mathcal{G}^S$ with zero mean and covariance matrix $L^{-1}$ (the inverse of the Laplacian matrix). This step is performed by employing the Cholesky factorization of the covariance matrix, $L^{-1} = \Omega^\top \Omega$. The simulated values are obtained in the following way
        $$ \left(Z_\mathcal{G}^S(v_1),\hdots, Z_\mathcal{G}^S(v_n)\right)^\top =  \Omega^\top Z_0, $$
        where $Z_0$ is an $n$-dimensional vector with uncorrelated standard Gaussian random variables.
        
   \item[] {\bf Step 2.} We extend the simulation on the vertices to locations on the whole graph. This is performed through linear interpolation (see Equation (\ref{interpolation})).

    %\item[] {\bf Step 3.} Finally, at each edge $e$, we add an independent Brownian bridge $Z_e(u)$. We use the following well-known Fourier expansion
        %\begin{equation}
        %    \label{spectral-bb}
        %    Z_e(u) = \sum_{k=1}^\infty \xi_k \frac{\sqrt{2(\overline{e}-\underline{e})}\sin(k\pi \varphi_e(u)/(\overline{e}-\underline{e}))}{k \pi}, \qquad u\in e,
        %\end{equation}
%where $\xi_1,\xi_2,\hdots$ are independent identically distributed standard Gaussian random variables.

    \item[] {\bf Step 3.}  {Finally, at each edge $e$, we add an independent Brownian bridge $B_e^S$, i.e., a Gaussian {random field} having a linear semi-variogram with slope $1/2$ and conditioned to $B_e^S(\underline{e})=B_e^S(\overline{e})=0$. To this end, we discretize $[\underline{e},\overline{e}]$ into a set of points $\underline{e} = \varepsilon_{0} < \varepsilon_{1} < \ldots < \varepsilon_{K_e-1} < \varepsilon_{K_e} = \overline{e}$ and use a sequential approach that takes advantage of the Markov property of the Brownian bridge \citep[p. 492]{Chiles2012}: for each $k \in [1,K_e-1]$, $B_e^S(\varepsilon_k)$ is a Gaussian random variable with mean equal to the ordinary kriging prediction from the leftmost and rightmost adjacent values $B_e^S(\varepsilon_{k-1})$ and $B_e^S(\overline{e})$, and variance equal to the corresponding ordinary kriging variance. This yields $$Z_e^S(u_k) = B_e^S(\varepsilon_{k}) = \left[\frac{1}{2} + \frac{\overline{e}-2\varepsilon_k+\varepsilon_{k-1}}{2(\overline{e}-\varepsilon_{k-1})}\right] B_e^S(\varepsilon_{k-1}) + \xi_k \sqrt{\frac{\overline{e}-\varepsilon_{k-1}}{4}-\frac{(\overline{e}-2\varepsilon_k+\varepsilon_{k-1})^2}{4(\overline{e}-\varepsilon_{k-1})}},$$ with $\xi_k$ a standard Gaussian random variable independent of $\xi_1, \ldots,\xi_{k-1}$, and $\{u_k = \varphi_e^{-1}(\varepsilon_k): k = 1,\ldots,K_e-1 \}$ is a set of points discretizing $e$.}

    \end{itemize}

\begin{remark}
Let us highlight some remarks on this methodology.
\begin{enumerate}
    \item  {$Z_\mathcal{G}^S$ is a Gaussian random field with zero mean and variogram given by the resistance metric, as per (\ref{metric}). Accordingly, it is an extension to the graph endowed with the resistance metric of the Brownian motion (Wiener-L\'evy process) on Euclidean spaces endowed with the Euclidean distance.}
    \item In real applications, the number of vertices is typically less than a few hundreds (e.g., $n=119$ and $n = 643$ for Eastbourne and Medellin traffic accident data, respectively, see \citealp{Moradi}), so the Cholesky approach is computationally tractable. Furthermore, although multiple independent copies of this auxiliary {random field} must be simulated (as described in the algorithms below), the Cholesky decomposition of $L^{-1}$ is performed only once. This is because the graph's structure, and consequently the matrix $L^{-1}$, remains constant throughout.
    
    %\item To implement Step 3,  we employ an appropriate truncation of (\ref{spectral-bb}).  {Alternativa: usar simulación por muestreador de Gibbs, ver \cite{Arroyo}, para simular un vector Gaussiano grande de covarianza conocida; evita la truncación que introduce un error. Mejor: simulación secuencial aprovechando la propiedad markoviana del movimiento Browniano, ver \cite[p. 492--493]{Chiles2012}}
    \item For each edge $e$, $Z_\mathcal{G}^S$ can be simulated with efficiency on a dense grid  {(i.e., $K_e$ may be very large)}. As a consequence, the methods introduced below will provide efficient algorithms for simulating Gaussian {random fields} on a large number of points on graphs.
\end{enumerate}
\end{remark}

\subsection{Algorithm 1 (Spectral Method)}

The first algorithm is based on the construction (\ref{eq1}). As shown in Proposition \ref{prop_spectral}, this approach allows us to obtain  {random fields} with completely monotone correlation structures. However, these {random fields} are clearly non-Gaussian, although their univariate distributions are standard Gaussian for any $u \in {\cal G}$ (Proposition \ref{prop_spectral2}). %\textcolor{red}{In particular, arguments in \cite{allard2020} prove that the distributions of pairs $(Y(u),Y(v))$ (with $u, v \in {\cal G}$) are mixtures of bivariate Gaussian distributions and have an isofactorial representation with the Hermite polynomials as the factors (see also \cite{matheron1982, matheron1989}).}
In order to obtain an approximately Gaussian {random field}, we consider an additive combination of several independent copies of (\ref{eq1}), as detailed in Algorithm \ref{algo:spectral},  {based on the central limit theorem \citep{lantu2002}}.

%\begin{equation}
    %\label{spectral-T}
    %Y^S(u) = \sqrt{\frac{2}{M}} \sum_{m=1}^M  {\sqrt{-\ln(V_m)}} \cos\left(W_m Z^{(m)}_{\mathcal{G}}(u) +  \Lambda_m\right), \qquad u\in\mathcal{G},
%\end{equation}
  %where $M$ is a large integer, $W_1,\hdots,W_M\sim F$,  {$V_1,\hdots,V_M\sim \text{Unif}(0,1)$} and $\Lambda_1,\hdots,\Lambda_M\sim \text{Unif}(0,2\pi)$ are independent, and $Z^{(1)}_{\mathcal{G}},\hdots,Z^{(M)}_{\mathcal{G}}$ are independent copies of the process $Z_{\mathcal{G}}^S$.

\begin{algorithm}
\caption{Spectral simulation. The simulated random variables $W_1,\hdots,W_M$, $V_1,\hdots,V_M$, and $\Lambda_1,\hdots,\Lambda_M$ are mutually independent and independent of $Z^{(1)}_{\mathcal{G}},\hdots,Z^{(M)}_{\mathcal{G}}$}
\begin{algorithmic}[1]
\Require A large integer $M$
\Require A spectral measure $F$
\Require $M$ independent copies $Z^{(1)}_{\mathcal{G}},\hdots,Z^{(M)}_{\mathcal{G}}$ of the auxiliary  {random field} $Z_{\cal G}^S$ simulated on a fine grid $\{u_p: p = 1,\ldots, P \}$ discretizing ${\cal G}$
\State Simulate $W_1,\hdots,W_M\sim F$
\State Simulate $V_1,\hdots,V_M\sim \text{Unif}(0,1)$
\State Simulate $\Lambda_1,\hdots,\Lambda_M\sim \text{Unif}(0,2\pi)$
\For{$p = 1$ to $P$}
\State Calculate $Y^S(u_p) = \sum_{m=1}^M \sqrt{\frac{-2\ln(V_m)}{M}} \cos\left(W_m Z^{(m)}_{\mathcal{G}}(u_p) +  \Lambda_m\right)$ 
\EndFor
\end{algorithmic}
\label{algo:spectral}
\end{algorithm}

\subsection{Algorithm 2 (Poisson Dilution)}

The second algorithm is based on the construction (\ref{random_token}). As shown in Proposition \ref{prop_dilution}, this approach allows us to obtain {random fields} with a variety of correlation structures. 
%Once selected the function $f$, the key step is to simulate the Poisson point process on the interval $I\subset \mathbb{R}$:
%\begin{enumerate}
    %\item Simulate $N \sim \text{Poisson}(|I|)$. 
    %\item If $N>0$, then generate $x_1,\hdots,x_N \sim \text{Unif}(I)$ independently. 
    %\item Return $X=\{x_1,\hdots,x_N\}$.
%\end{enumerate}
Once selected the dilution function $f$, a direct application of (\ref{random_token}) requires defining a sufficiently large interval $I$ to simulate the Poisson point process $X \cap I$, which yields Algorithm \ref{algo:dilution} below.

\begin{algorithm}
\caption{Poisson dilution. The simulated random variables ($N, x_1, \ldots, x_N, \epsilon_1, \ldots, \epsilon_N)$ are mutually independent and independent of $Z_{\cal G}^S$}
\begin{algorithmic}[1]
\Require A bounded interval $I \subset \mathbb{R}$
\Require A dilution function $f$
\Require An auxiliary {random field} $Z_{\cal G}^S$ simulated on a fine grid $\{u_p: p = 1,\ldots, P \}$ discretizing ${\cal G}$
\State Simulate $N \sim \text{Poisson}(|I|)$
\If{$N>0$}
\State Simulate $x_1,\ldots,x_N \sim \text{Unif}(I)$
\State Simulate $\epsilon_1,\ldots,\epsilon_N \sim \text{Rademacher}$
\For{$p = 1$ to $P$}
\State Calculate $Y^S(u_p) = \sum_{n=0}^{N} \epsilon_n \, f(Z_{\cal G}^S(u_p) - x_n)$ 
\EndFor
\Else 
\State Deliver $Y^S(u_p) = 0$ for $p = 1,\ldots,P$
\EndIf
\end{algorithmic}
\label{algo:dilution}
\end{algorithm}

Because $I$ is bounded, Algorithm \ref{algo:dilution} is approximate, i.e., it does not exactly reproduce the target covariance structure, unless the dilution function $f$ is compactly supported, as in the first entry of Table \ref{tab:cov2}. 
In the same manner as for the spectral simulation, a central limit approximation is required to obtain a {random field} whose finite-dimensional distributions are close to multivariate normal.

\subsection{Algorithm 3 (Dilution of a Random Germ)}

Algorithm \ref{algo:importancesampling2} is a direct application of the third proposal in Section \ref{thirdconstruction}, combined with a central limit approximation to get a random {field} with finite-dimensional distributions close to multivariate normal. In practice, the importance sampling density $\varpi$ should be chosen in such a way that $x$ has a high probability to be in, or close to, the target simulation domain, and a low probability to be far from this domain.

\begin{algorithm}
\caption{Dilution of a random germ. The simulated random variables ($x_1,\ldots,x_M,\epsilon_1,\ldots,\epsilon_M)$ are mutually independent and independent of $Z^{(1)}_{\mathcal{G}},\hdots,Z^{(M)}_{\mathcal{G}}$}
\begin{algorithmic}[1]
\Require An importance sampling density $\varpi$ supported in $\mathbb{R}$
\Require A dilution function $f$
\Require A large integer $M$
\Require $M$ independent copies $Z^{(1)}_{\mathcal{G}},\hdots,Z^{(M)}_{\mathcal{G}}$ of the auxiliary {random field} $Z_{\cal G}^S$ simulated on a fine grid $\{u_p: p = 1,\ldots, P \}$ discretizing ${\cal G}$
\State Simulate $x_1,\ldots,x_M \sim \varpi$
\State Simulate $\epsilon_1,\ldots,\epsilon_M \sim \text{Rademacher}$
\For{$p = 1$ to $P$}
\State Calculate $Y^S(u_p) = \sum_{m=1}^M \frac{\epsilon_m}{\sqrt{M \varpi(x_m)}}  f(Z_{\cal G}^S(u_p) - x_m)$ 
\EndFor
\end{algorithmic}
\label{algo:importancesampling2}
\end{algorithm}

\section{Numerical Example}
\label{ejemplos}

% {Imagino lo siguiente. (1) Simulación no condicional: probar el método espectral, con o sin el término $-\ln(V)$, y comparar las velocidades de convergencia a un proceso Gaussiano, en función de $M$ (estilo https://dl.acm.org/doi/abs/10.1145/3421316). Para el método de dilución, comparar la versión truncada con la versión de muestreo de importancia con diferentes distribuciones de importancia $\varpi$ (binomial negativa desfasada de $1$ unidad, zeta, etc.) y diferentes $M$ para el límite central. (2) También se puede hablar de condicionar la simulación a datos existentes, vía kriging.}

\subsection{Implementation Parameters}

We apply our algorithms to the network of the University of Chicago neighborhood (see \citealp{baddeley2021analysing}),
which consists of 338 vertices and 503 edges. 

For Algorithm \ref{algo:spectral}, we consider the first entry of Table \ref{tab:cov} (an exponential correlation function), whereas for Algorithms \ref{algo:dilution} and \ref{algo:importancesampling2}, we consider the second entry of Table \ref{tab:cov2} (a Cauchy correlation function). The scale parameter $a$ is set to $0.2$ and $M$ is set to {$1000$} for the central limit approximation. In Algorithm \ref{algo:dilution}, we use the interval {$I=[-50,50]$}, while Algorithm \ref{algo:importancesampling2} is implemented with a standard Cauchy importance sampling density, favoring the occurrence of a wide range of values, attributed to its heavy tails.

\begin{figure}
    \centering
    \includegraphics[scale=0.15]{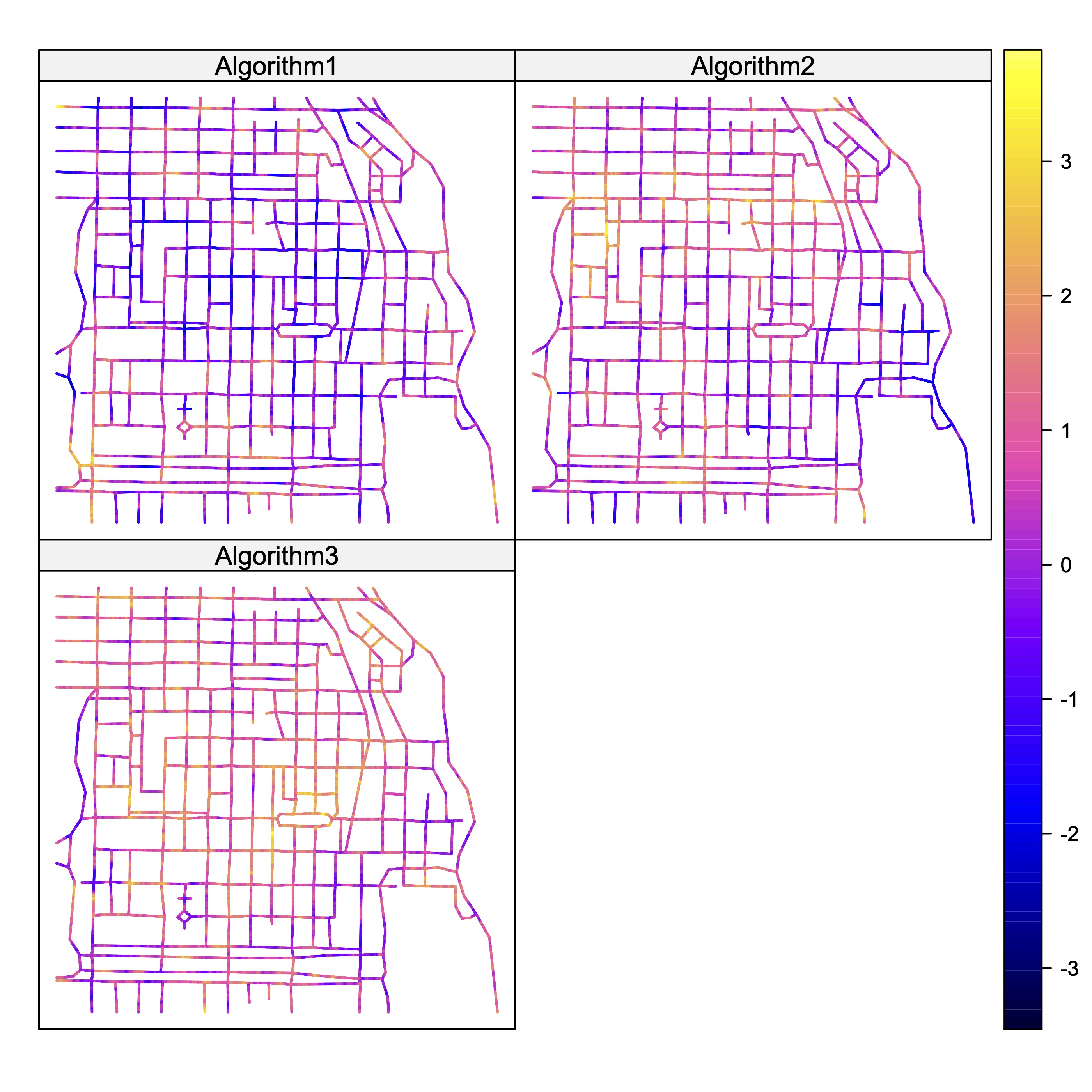}
    \caption{Realizations of random fields in a neighborhood of the University of Chicago, obtained from the three proposed algorithms, considering 200 equispaced points on each of the 503 edges (totaling 100,600 locations). For Algorithm \ref{algo:spectral}, we consider an exponential correlation function, whereas for Algorithms \ref{algo:dilution} and \ref{algo:importancesampling2}, we consider a Cauchy correlation function. In each realization, the scale parameter is $a=0.2$ and the central limit approximation is implemented with {$M=1000$}.}
    \label{fig:chicago1}
\end{figure}

As an illustration, Figure \ref{fig:chicago1} shows realizations over 100,600 locations (200 per edge) obtained by each algorithm.

\subsection{Execution Times}

We assess the numerical complexity of the algorithms with respect to the number of target points over the network. Since this network contains 503 edges, we consider $503 \times 2^k$ points over $\mathcal{G}$, with $k=5,6,\hdots,10$, i.e., $2^k$ represents the number of points considered for discretizing each edge. In this experiment, we used an Apple M1 Pro chip (10-core CPU) with 16GB of memory. The results (Table \ref{tab:tiempos}) indicate that Algorithms \ref{algo:spectral} and \ref{algo:importancesampling2} exhibit comparable execution times. The execution times of Algorithm \ref{algo:dilution} are greater because the interval {$I=[-50,50]$}, which is necessary for a good coverage, implies the generation of a Poisson process with a large number $N$ of points (the expectated value of $N$ is {$|I|=100$}). Therefore, constructing the random field is computationally more demanding than in the other proposed algorithms.

\begin{table}
    \centering
    \caption{Computing times (in seconds) to simulate a  random field, as a function of the number of target points over the network, for each algorithm, with {$M=1000$} in the central limit approximation.}
    
    \begin{tabular}{ccccccc} \hline \hline
       & \multicolumn{6}{c}{Number of target points} \\ \cline{2-7}
                       &16,096 &32,192 & 64,384 &128,768 &257,536 & 515,072\\ \hline
       Algorithm \ref{algo:spectral}  &1.837&2.762&5.431&10.61&21.49&48.86 \\
       Algorithm \ref{algo:dilution}  &29.75&66.36&134.8&274.1&604.1&1263\\
       Algorithm \ref{algo:importancesampling2} &1.408&2.396&4.431&8.701&17.91&42.49 \\ \hline \hline
    \end{tabular}
    \label{tab:tiempos}
\end{table}

\subsection{Correlation Structure Reproduction}

\begin{figure}
    \centering
    \includegraphics[scale=0.053]{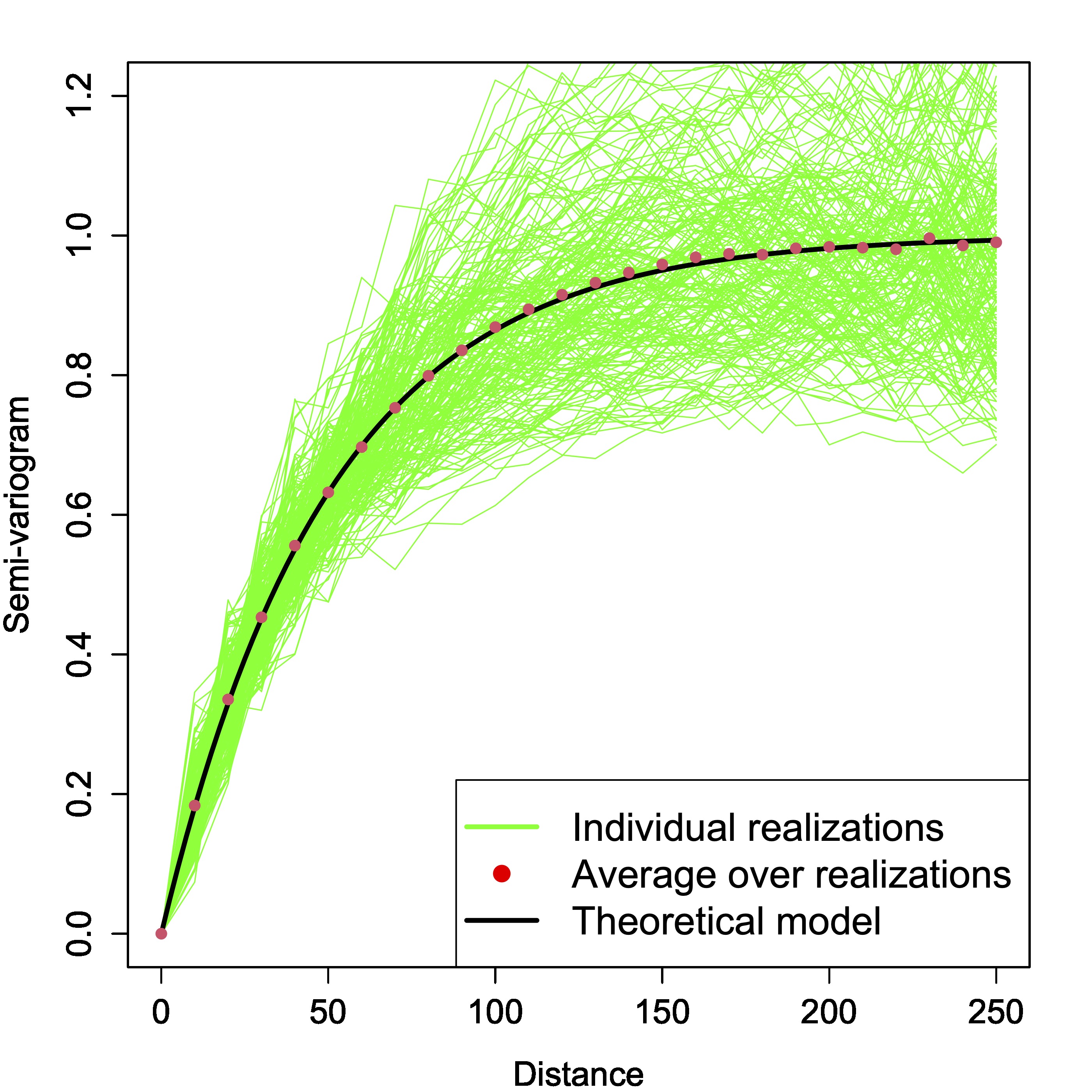}
    \includegraphics[scale=0.053]{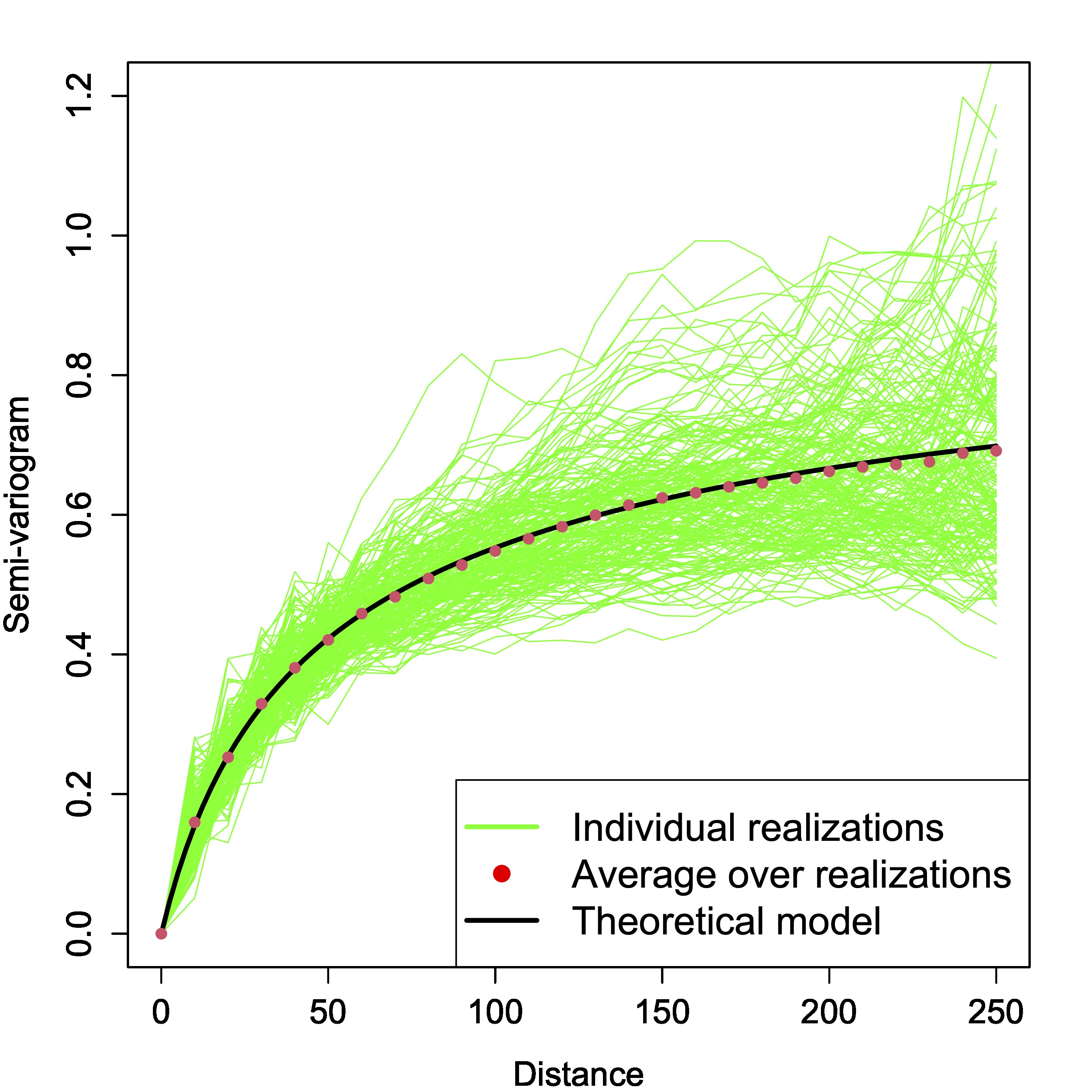}
    \includegraphics[scale=0.053]{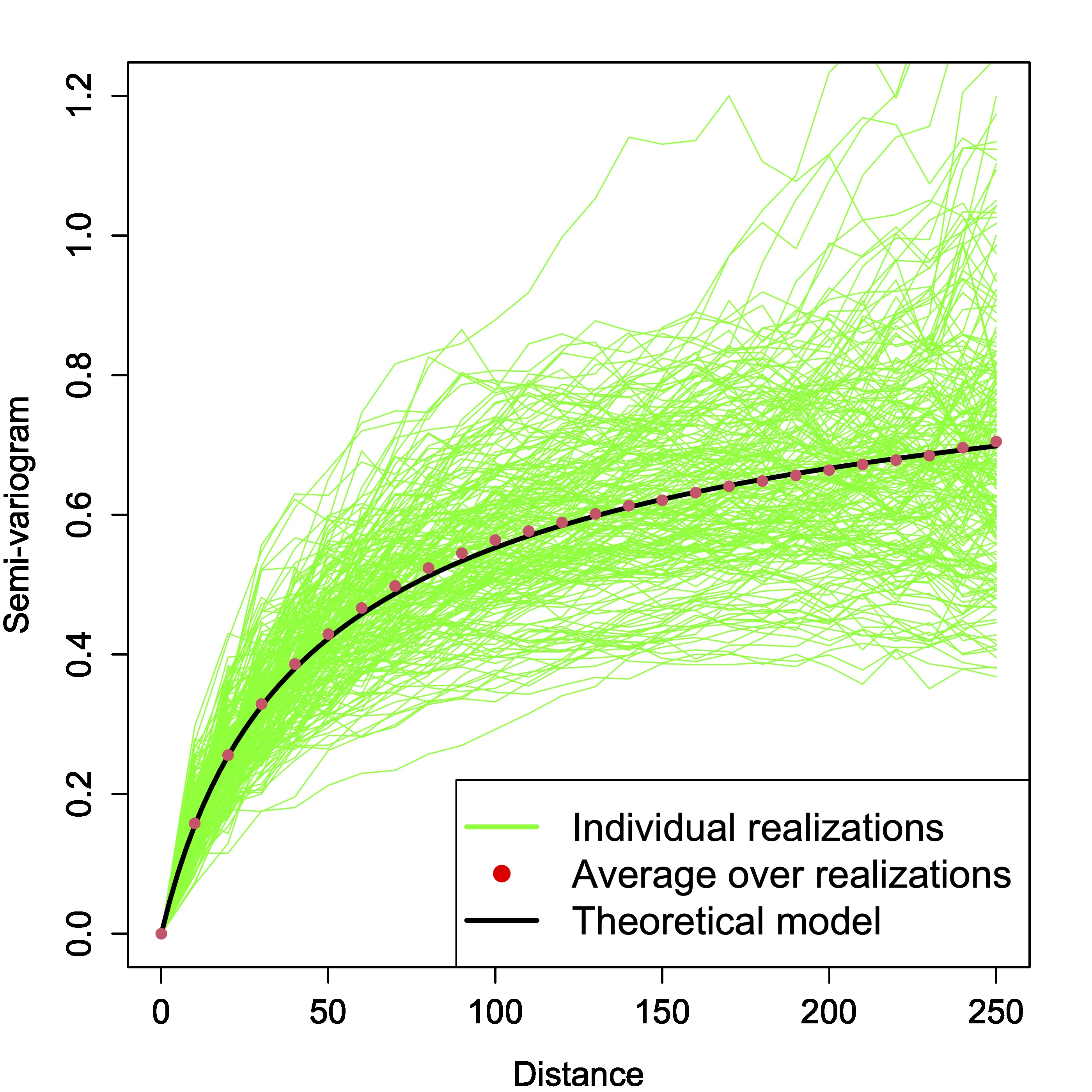}
    \caption{Experimental (green lines) and theoretical (black lines) semi-variograms for {200} realizations of a random field simulated on 2 points per edge ({1006} points in total), as a function of the resistance metric. The red points are the average of the experimental semi-variograms. From left to right, we consider Algorithms \ref{algo:spectral} to \ref{algo:importancesampling2}, respectively.}
    \label{fig:variograms1}
\end{figure}

To investigate the reproduction of the target spatial dependency structure, for each algorithm, we construct $200$ realizations at a set of locations discretizing the graph (2 points per edge, i.e., 1006 points in total), with the same random field models and parameters as in the previous subsections, and compute their empirical semi-variograms. It is seen (Figure \ref{fig:variograms1}) that, on average, the theoretical semi-variogram is reproduced without any bias, a conclusion that is confirmed by statistical testing (Table \ref{tab:t_test_variograms}) \citep{emery2008}.

%\textcolor{green}{Podría hacer un test T-Student sobre los variogramas para algunas distancias de referencia [10, 50, 100, 250] y agregar los resultados en una tabla. Esto dará un poco más de ``carne'' a esta sección.}

\begin{table}
    \centering
    \begin{tabular}{ccccccc} \hline \hline
       & \multicolumn{4}{c}{Absolute value of T-statistics} \\ \cline{2-7}
                       &Lag 10 &Lag 50 &Lag 100 & Lag 150 & Lag 200 & Lag 250 \\ \hline
       Algorithm \ref{algo:spectral} & 0.661 & 0.010 & 0.533 & 0.899 & 0.232 & 0.249\\
       Algorithm \ref{algo:dilution} & 1.602 & 0.627 & 1.052 & 0.342 & 0.619 & 0.625\\
       Algorithm \ref{algo:importancesampling2} & 0.991 & 1.131  & 1.426 & 0.147 & 0.275 & 0.461\\ \hline \hline
    \end{tabular}
    \caption{ {Student test on experimental semi-variograms at five lag distances ($d_R(u,v) = 10$, $50$, $100$, $150$,  $200$ and $250$), for $200$ realizations of the simulated random field. The critical value at a $0.05$ level of significance is $1.972$. The null hypothesis that the average experimental semi-variogram matches the theoretical semi-variogram is accepted in all the cases.} }
    \label{tab:t_test_variograms}
\end{table}

\begin{figure}
   \centering
    \includegraphics[scale=0.053]{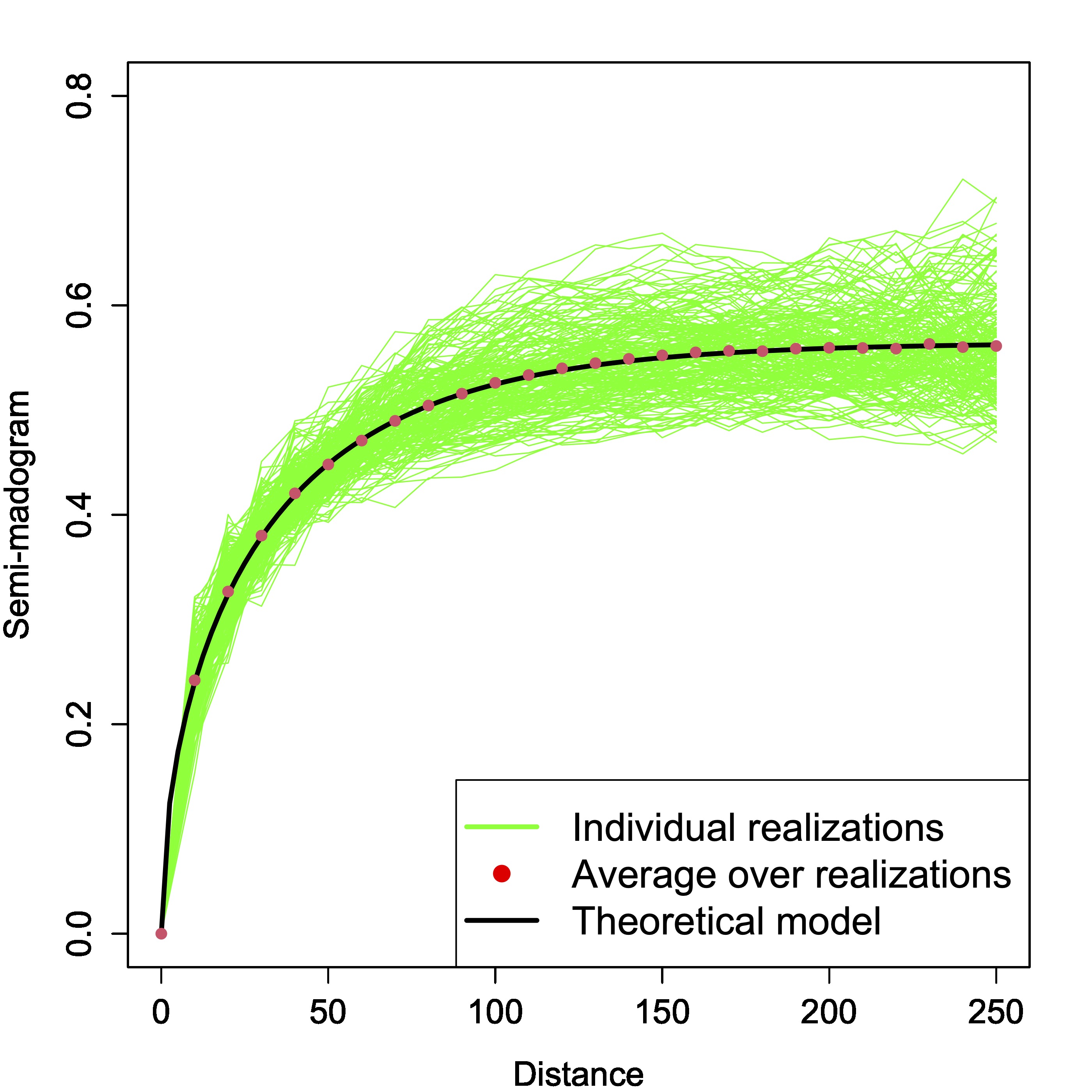}
    \includegraphics[scale=0.053]{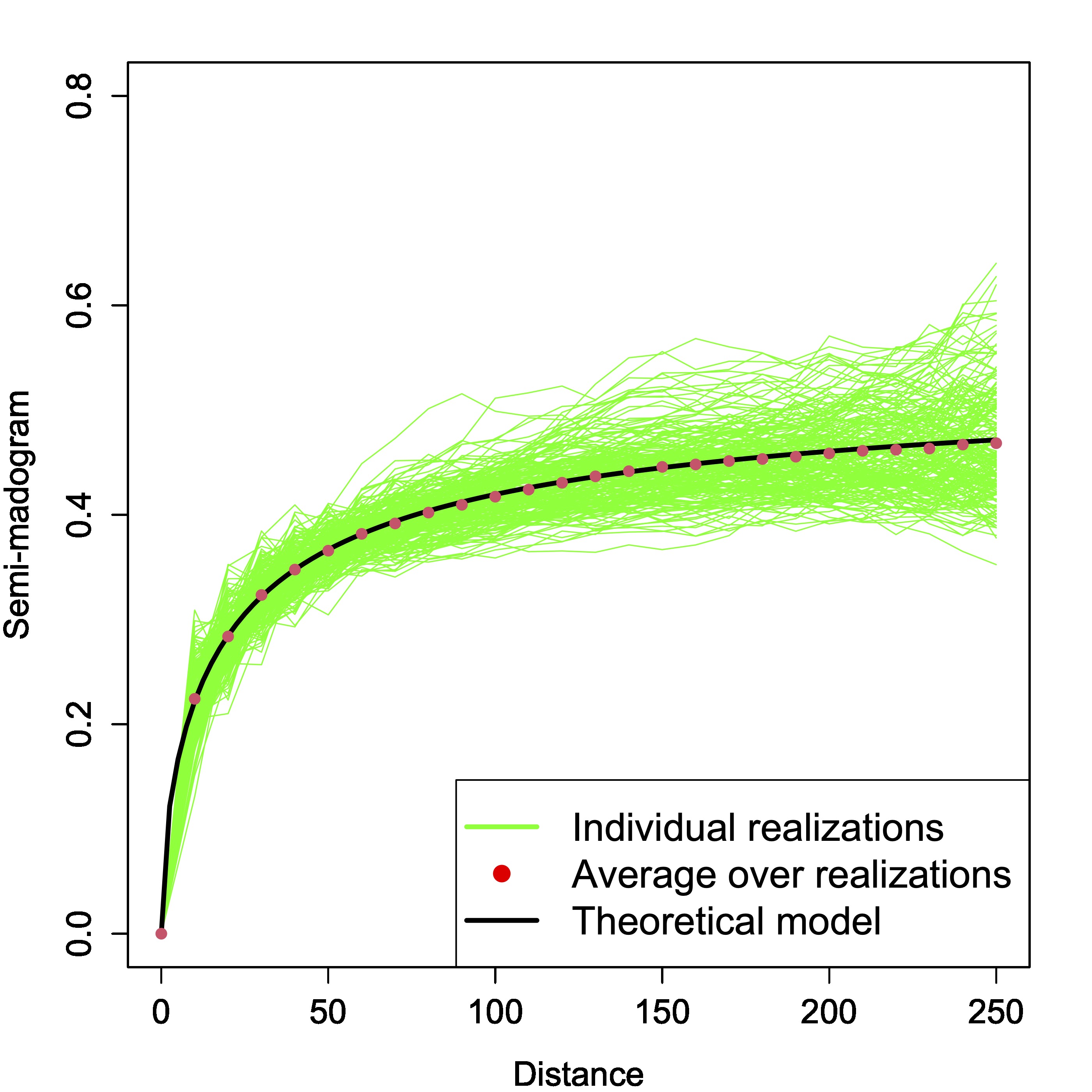}
    \includegraphics[scale=0.053]{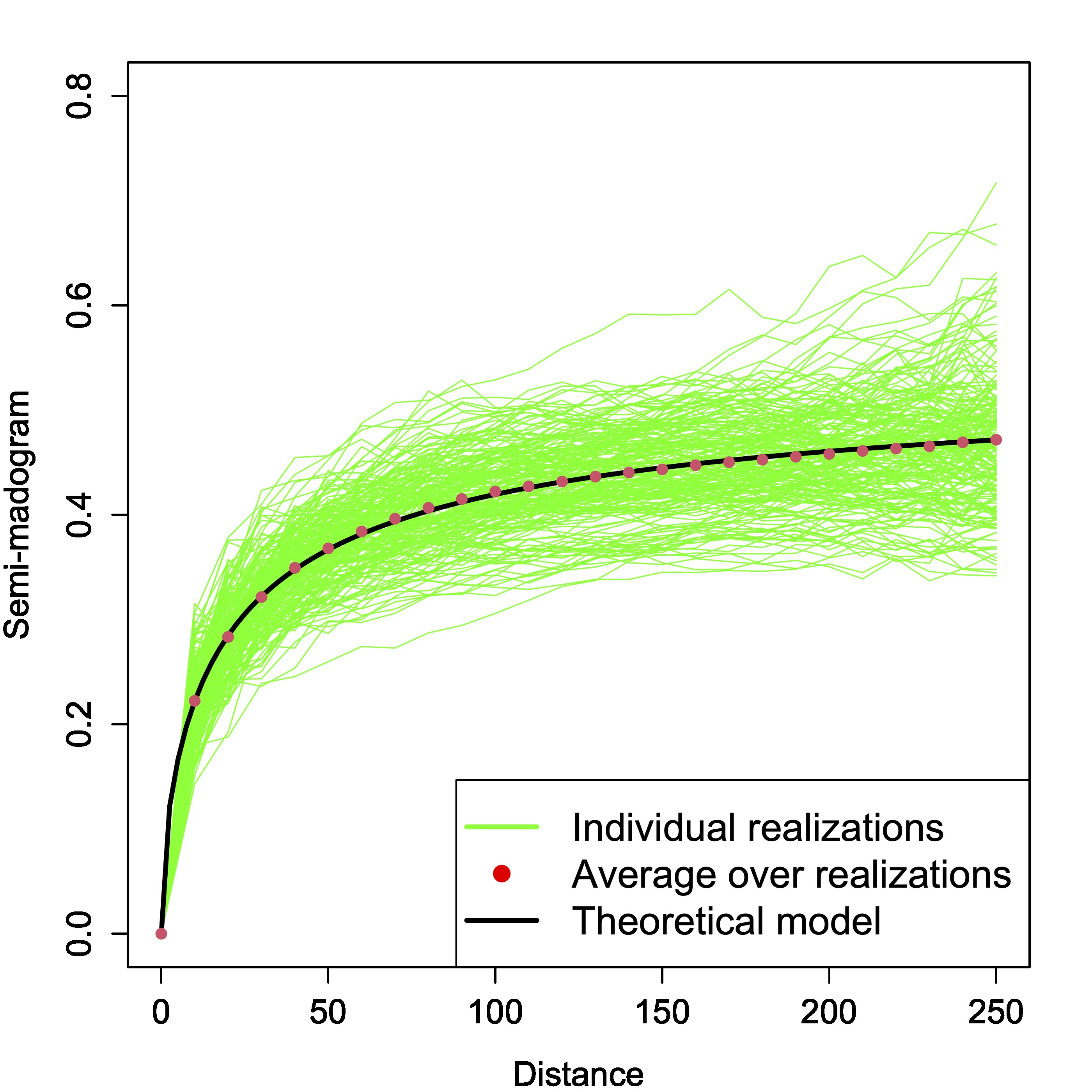}
    \caption{{Experimental (green lines) and theoretical (black lines) semi-madograms for {200} realizations of a random field simulated on {2} points per edge ({1006} points in total), as a function of the resistance metric. The red points are the average of the experimental semi-madograms. From left to right, we consider Algorithms \ref{algo:spectral} to \ref{algo:importancesampling2}, respectively.}}
    \label{fig:madograms}
\end{figure}

{The same conclusions prevail when looking at the semi-madograms of the simulated random fields (Figure \ref{fig:madograms} and Table \ref{tab:t_test_madograms}): on average over the 200 realizations, the experimental semi-madograms match the expected model and no bias is perceptible. Recall that, for a Gaussian random field, the semi-madogram $\gamma_1$ is proportional to the square root of the semi-variogram $\gamma_2$ \citep{emery2005}:  
$$\gamma_1 = \sqrt{\frac{\gamma_2}{\pi}} = \sqrt{\frac{C_Y(0)-C_Y}{\pi}},$$
with $\gamma_1(d_R(u,v)) = \frac{1}{2} \mathbb{E}\left[\mid Y(u)-Y(v)\mid \right]$ and $\gamma_2(d_R(u,v)) = \frac{1}{2} \mathbb{E}\big[\left(Y(u)-Y(v) \right)^2\big]$ for $u,v \in {\cal G}$.} 

In passing, the fact that the semi-madogram is well reproduced indicate that the bivariate distributions of the simulated random fields are (close to) bigaussian, i.e., the central limit approximation is accurate.

\begin{table}
    \centering
    \caption{{Student test on experimental semi-madograms at five lag distances ($d_R(u,v) = 10$, $50$, $100$, $150$,  $200$ and $250$), for $200$ realizations of the simulated random field. The critical value at a $0.05$ level of significance is $1.972$. The null hypothesis that the average experimental semi-madogram matches the theoretical semi-madogram is accepted in all the cases.} }

    \begin{tabular}{ccccccc} \hline \hline
       & \multicolumn{4}{c}{Absolute value of T-statistics} \\ \cline{2-7}
                       &Lag 10 &Lag 50 &Lag 100 & Lag 150 & Lag 200 & Lag 250 \\ \hline
       Algorithm \ref{algo:spectral} & 0.757 & 0.316 & 0.570 & 0.866 &  0.193 & 0.319\\
       Algorithm \ref{algo:dilution} & 1.030  & 0.819 & 1.307 & 0.336 & 0.781 & 0.831\\
       Algorithm \ref{algo:importancesampling2} & 0.130  & 0.474  & 0.932 & 0.552 & 0.747 & 0.033\\ \hline \hline
    \end{tabular}
    \label{tab:t_test_madograms}
\end{table}

\subsection{Assessment of Central Limit Approximation}

To end the case study, we explore the accuracy of the multivariate normal approximation in terms of the number $M$ of independent copies in the central limit approximation.  {Both analytical and numerical approaches have been developed in the literature to this end. However, the former (e.g., \citealp{lantu1994}) are impracticable in our context, due to the difficulty in deriving an analytical expression of the finite-dimensional distributions of the simulated random fields, or a Berry-Esseen bound of the error between the true and simulated distributions. Here, we opt for a numerical validation based on statistical testing and inspired from \cite{Arroyo2}, as detailed hereinafter. }

\begin{figure}
    \centering
    \includegraphics[scale=0.07]{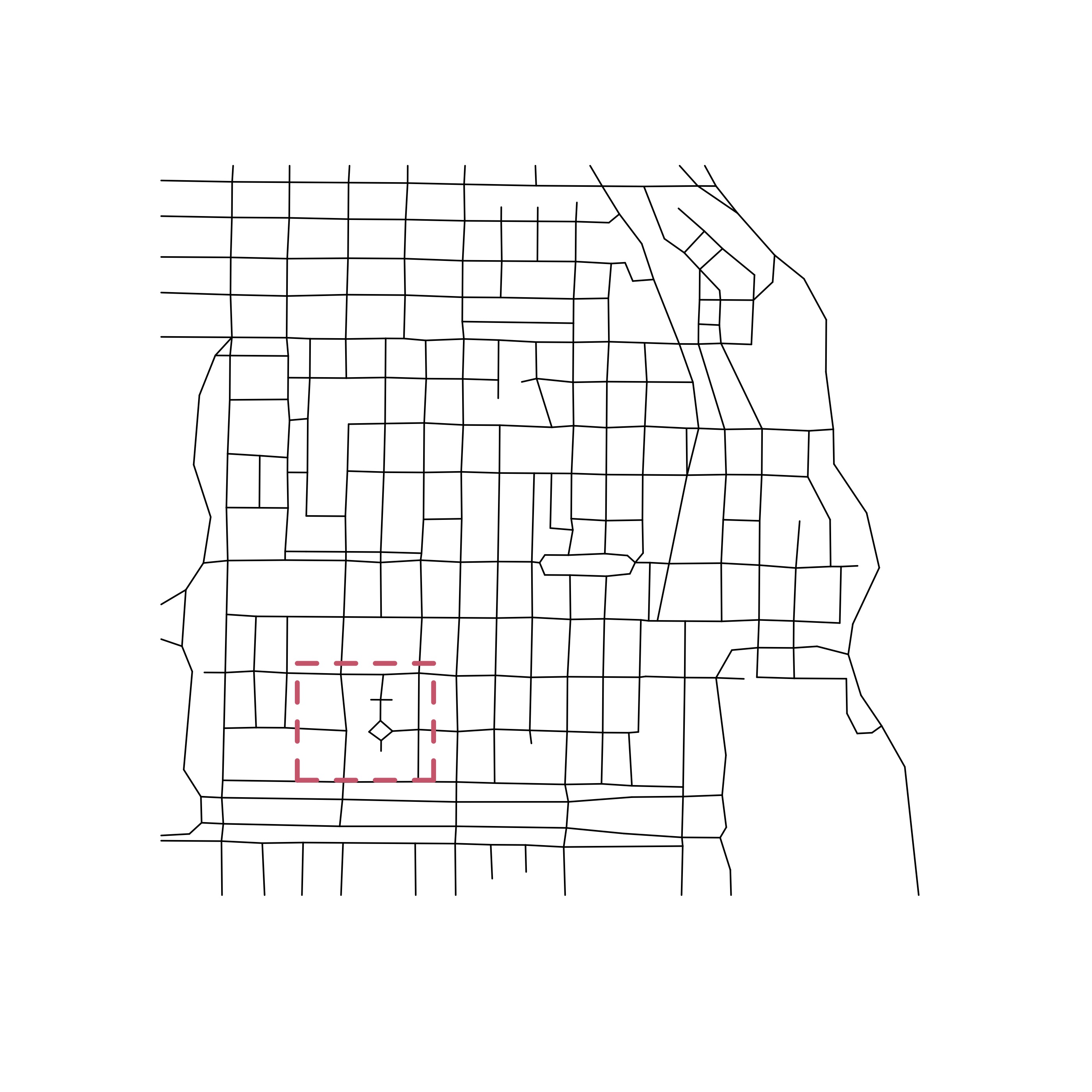}\includegraphics[scale=0.07]{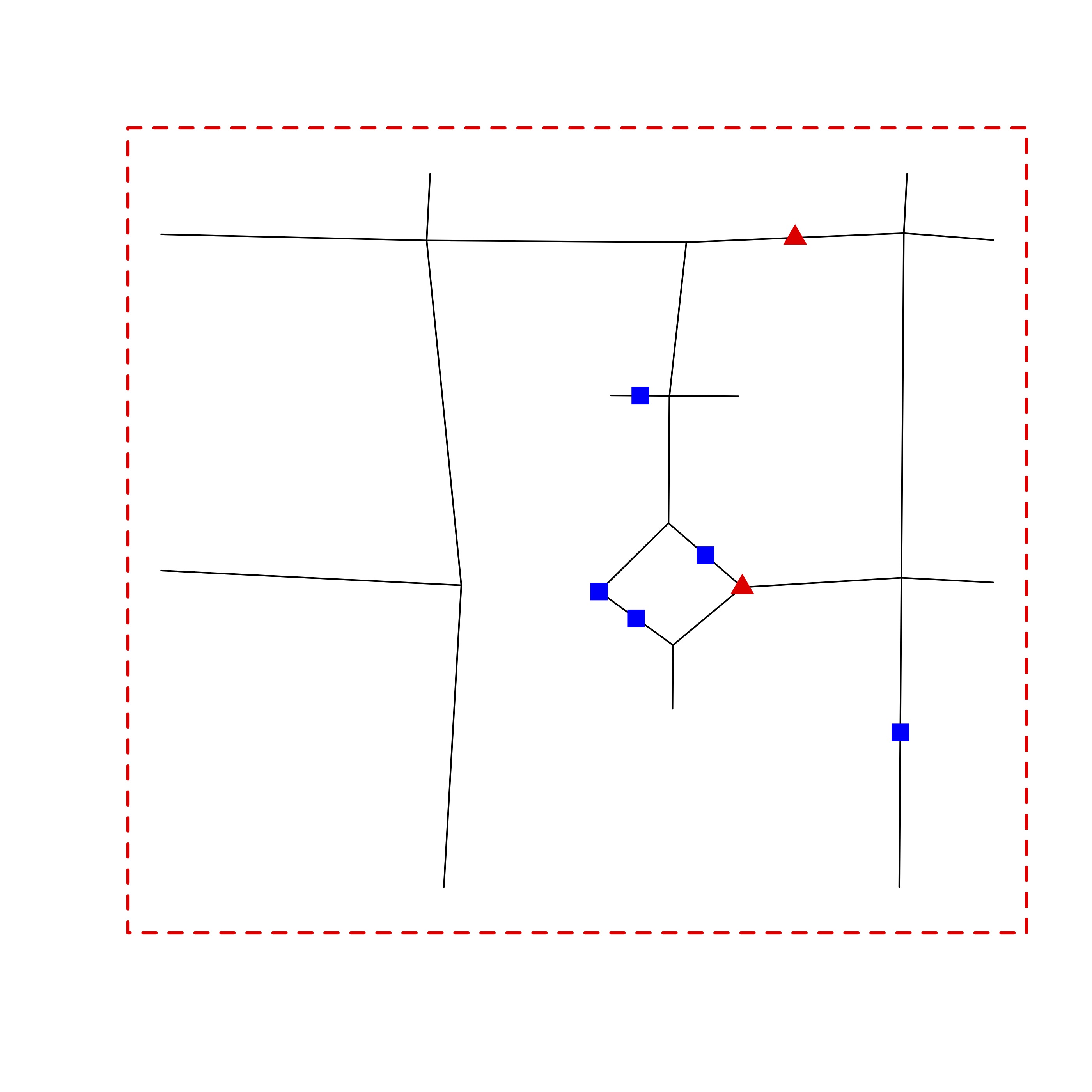}
    \caption{(Left) The dashed red square shows the zone of the University of Chicago neighborhood employed in our experiment. (Right) Locations within this zone targeted for simulation.}
    \label{fig:simplenet}
\end{figure}

Given $n$ locations on the graph, we study the distribution of a linear combination of the values simulated at these locations. To simplify the analysis, {let us consider a subset of the University of Chicago neighborhood. We consider two scenarios: $n=2$ (red triangles in Figure \ref{fig:simplenet}) and $n=5$ (blue squares in Figure \ref{fig:simplenet}) fixed locations. Thus, in each scenario, we consider the zero-mean linear combination $\lambda_1 Y(u_1)+\cdots+\lambda_nY(u_n)$, where the weights $\lambda_1,\hdots,\lambda_n$ are simulated from a uniform distribution on $[-10,10]$. These weights are fixed across this experiment, therefore one can calculate the variance of the linear combination once the covariance function of $Y(\cdot)$ is known.}

\begin{figure}
    \centering
        \includegraphics[scale=0.1]{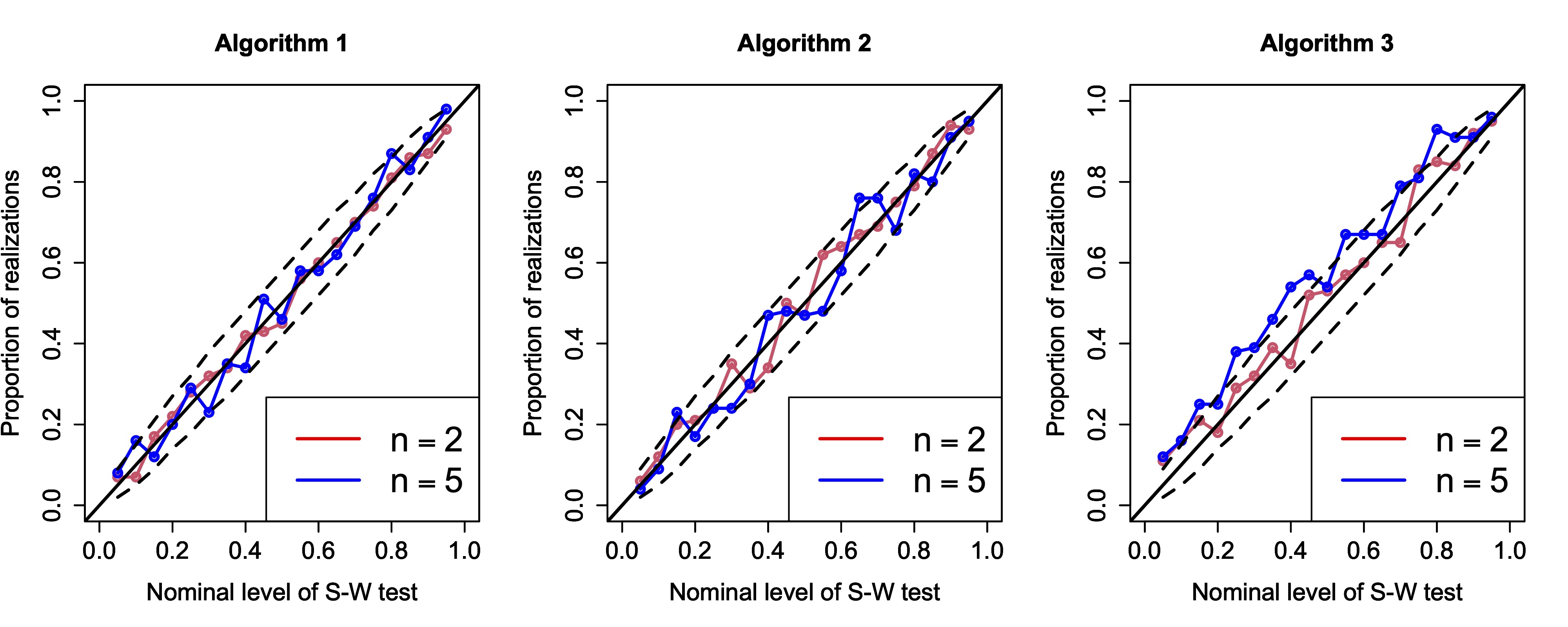}
    \includegraphics[scale=0.1]{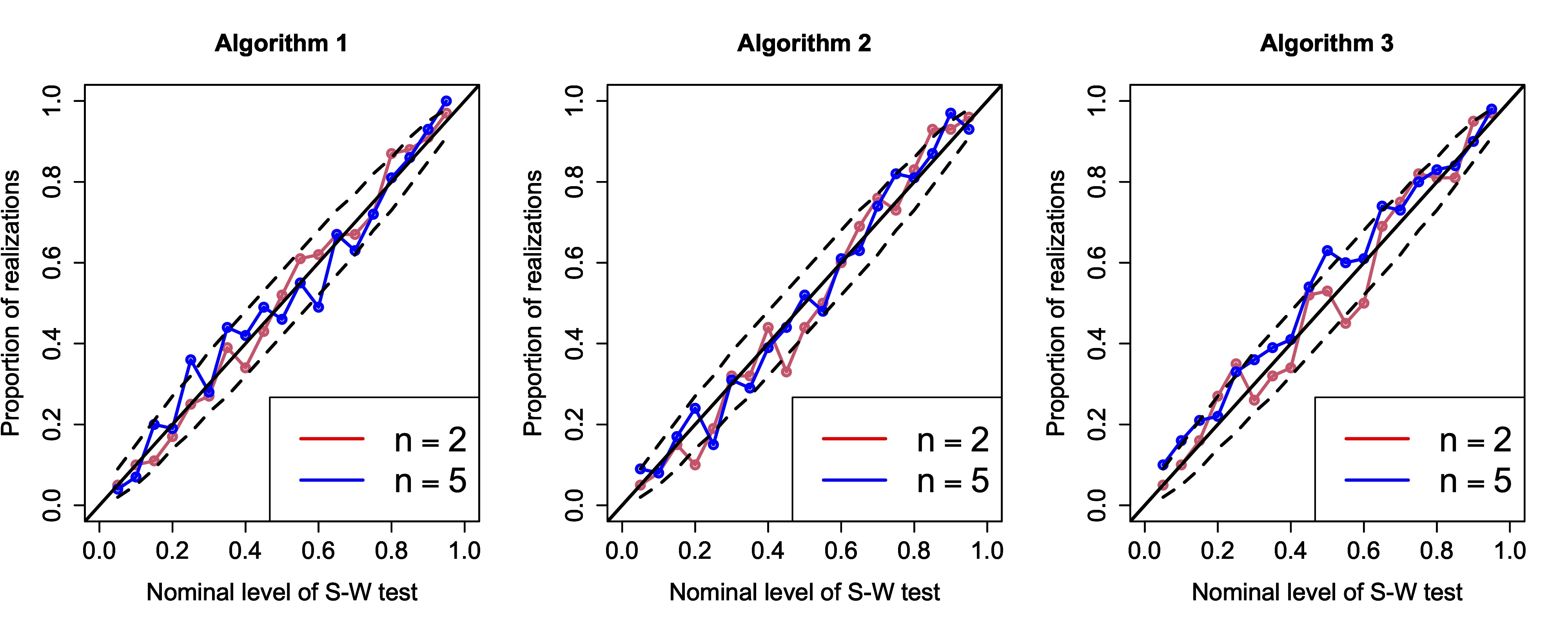}
    \includegraphics[scale=0.1]{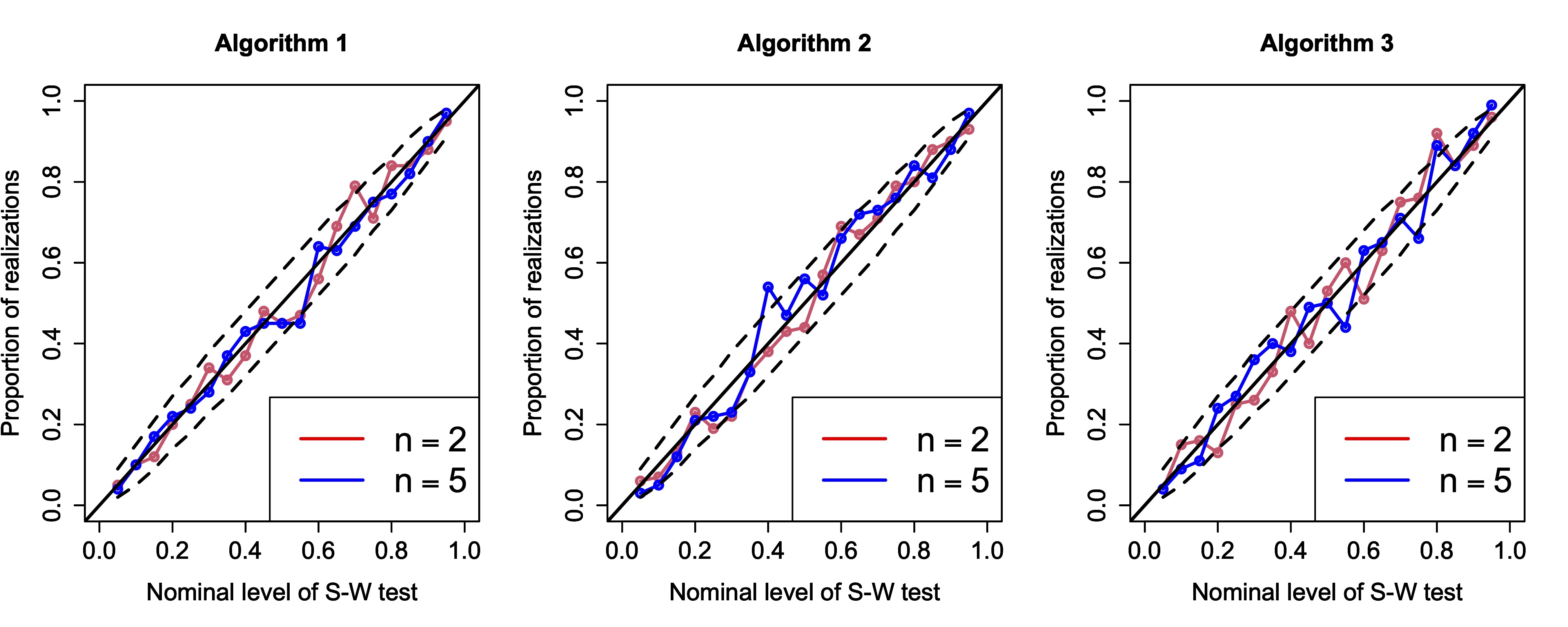}
    \caption{Probability-probability plots showing the proportion of rejected Shapiro-Wilk tests on linear combinations of simulated observations versus the nominal significance level of the test.  {We consider, from top to bottom,  $M=50$, $100$  and $500$.} In each panel, the black solid line is the identity line, whereas the dashed black lines represent the 90\% confidence bounds on the observed proportions of rejections.   }
    \label{fig:pp-plot}
\end{figure}

To assess multivariate normality, we proceed as follows. We simulate 100 independent linear combinations, as detailed above, and use the Shapiro-Wilk goodness-of-fit test, with a certain significance level $\alpha\in (0,1)$, to check the normal assumption. By repeating this procedure 100 times, we can calculate the proportion of rejections and compare this value with the nominal level $\alpha$. Figure \ref{fig:pp-plot} shows this comparison for several values of $\alpha$, and considering each algorithm with {$M=50$, $100$ and $500$}. We also report a $90\%$ confidence interval for the observed proportions of rejected tests, which is constructed by taking into account that, for each $\alpha$, the number of rejections follows a binomial distribution with size 100 and success probability $\alpha$.

Both algorithms \ref{algo:spectral} and \ref{algo:dilution} offer robust Gaussian approximations even with $M=50$, as the proportion of rejections almost always lies in-between the 90\% confidence bounds derived from the binomial distribution. Conversely, Algorithm \ref{algo:importancesampling2} requires raising $M$ to achieve a more accurate performance, as the linear combinations with $n=5$ appear to significantly deviate from normality. This experiment reveals that, for the network under consideration, $M$ of the order of a few hundreds is sufficient to achieve a good multivariate-Gaussian approximation across all algorithms.

{This kind of goodness-of-fit testing is deemed useful to determine a suitable value for $M$, especially for Algorithms \ref{algo:dilution} and \ref{algo:importancesampling2} where the quality of the multivariate-Gaussian approximation may strongly depend on the support of the dilution function and on the network geometry; in contrast, for Algorithm \ref{algo:spectral}, the univariate distribution of the simulated random field is always Gaussian (Proposition \ref{prop_spectral2}), thus multivariate normality can be reached more easily.}

\section{Discussion and Conclusions}
\label{discussion}

The algorithms presented in this study allows the simulation of Gaussian {random fields} on graphs with Euclidean edges, being isotropic with respect to the resistance metric. In developing these algorithms, we have introduced interesting stochastic representations that offer an alternative path to embeddings in Hilbert spaces for obtaining valid classes of covariance functions on these complex domains.

In particular, Algorithm \ref{algo:spectral}, which allows us to simulate {random fields} with completely monotone covariance functions, represents a graph-like version of classical spectral algorithms from geostatistics that have been applied in past literature to simulate random fields indexed by either Euclidean \citep{Shinozuka, mantoglou1982turning} or spherical \citep{alegria2020turning} coordinates. Algorithms \ref{algo:dilution} and \ref{algo:importancesampling2}, which adapt the dilution principle detailed in \cite{matern} and \cite{lantu2002}, allow us to achieve a wide range of covariance functions; however, Algorithm \ref{algo:dilution} is {approximate, unless the dilution function is compactly supported}, whereas Algorithm \ref{algo:importancesampling2}, which employs an importance sampling strategy, reproduces the target covariance structure exactly {irrespective of the dilution function under consideration. Note that the latter algorithm can be adapted to random fields defined in Euclidean domains as an alternative to the Poisson dilution method when the dilution function is not compactly supported, avoiding biases in the reproduction of the covariance structure due to the restriction of the Poisson process to a bounded interval.}

We demonstrated that the proposed algorithms reproduce the target correlation structures. {For all of them, the numerical complexity is, up to a fixed computational cost, proportional to the number of locations targeted for simulation on the graph, which ensures their applicability to large-scale problems. For instance, with standard computer resources,} Algorithms \ref{algo:spectral} and \ref{algo:importancesampling2} generate realizations across half a million points on the graph in about {$45$ seconds}, {whereas classical one-size-fits-all algorithms, such as the sequential Gaussian and Cholesky decomposition algorithms \citep{Chiles2012}, are out of reach}. The simulation of an auxiliary {random field} related to the resistance metric limits the speed slightly. However, the algorithms remain notably fast {compared to existing alternatives}. Statistical testing indicates that combining a few hundred independent copies yields satisfactory Gaussian approximations. 

The proposed stochastic representations also hold promise for developing covariance novel models and simulation algorithms in future investigations, particularly when addressing vector-valued {random fields}, for which the covariance function is matrix-valued, and spatio-temporal {random fields}.

{Our findings also provide valuable insights for simulating point pattern data on networks, specifically focusing on log-Gaussian Cox processes, where the logarithm of the intensity function is governed by a Gaussian random field \citep{cox1955some}; see \cite{baddeley2021analysing} for an inventory of potential applications in the social sciences, life sciences, physical sciences and engineering.}

\section*{Acknowledgments} 
This work was partially funded by the National Agency for Research and Development of Chile [grants ANID Fondecyt 1210050 (A. Alegr\'ia and X. Emery) and ANID PIA AFB230001 (X. Emery)]. Alfredo Alegr\'ia also acknowledges the support of Universidad Técnica Federico Santa María for grant Proyectos Internos USM 2023 PI$_{-}$LIR$_{-}$23$_{-}$11.

\bibliographystyle{apalike}
\bibliography{mybib}

\end{document}